\documentclass[10pt]{article}
\usepackage{float}
\usepackage[symbol]{footmisc}
\usepackage{graphicx}
\usepackage[superscript,sort]{cite}
\usepackage{array}
\usepackage{bm}
\usepackage{appendix}
\newcolumntype{x}[1]{%
>{\centering\hspace{0pt}}p{#1}}%
\usepackage[normalem]{ulem} 
\usepackage[colorlinks=true,linkcolor=black,citecolor=black,urlcolor=black,bookmarks=true,breaklinks=true]{hyperref}
\usepackage[usenames]{color}
\usepackage{amsmath,amssymb}

\usepackage{amsthm,amscd,amsxtra,amsfonts,amsmath,amssymb,multirow}
\usepackage{wrapfig}
\usepackage[footnotesize]{caption}
\usepackage{subcaption}
\usepackage[tiny,compact]{titlesec}
\usepackage{threeparttable} 

\usepackage{helvet}

\usepackage{booktabs}
\usepackage[table]{xcolor}
\usepackage{xcolor,colortbl}
\usepackage{placeins}
\usepackage{tikz}
\usetikzlibrary{shapes,arrows}
\usepackage[T1]{fontenc}
\usepackage[linesnumbered,ruled]{algorithm2e} 
\definecolor{MyPurple}{rgb}{1,0,1}

\newcommand{\beq}[1]{\begin{equation} \label{#1}}
\newcommand{\eeq}{\end{equation}}
\newcommand{\barray}{\begin{array}{ll}}
\newcommand{\earray}{\end{array}}

\usepackage{longtable} 
\usepackage{threeparttable} 
\usepackage{booktabs} 

\usepackage{tikz-cd}
\def\be{\mathbf{e}}
\def\bn{\mathbf{n}}
\def\bv{\mathbf{v}}

\def\bu{\mathbf{u}}
\def\ba{\mathbf{a}}
\def\bb{\mathbf{b}}
\def\bc{\mathbf{c}}
\def\bh{\mathbf{h}}

\def\bt{\mathbf{t}}
\def\br{\mathbf{r}}
\def\im{\mathrm{im}\;}
\makeatother
\tikzset{
	subseteq/.style={
		draw=none,
		edge node={node [sloped, allow upside down, auto=false]{$\subseteq$}}},
	Subseteq/.style={
		draw=none,
		every to/.append style={
			edge node={node [sloped, allow upside down, auto=false]{$\subseteq$}}}
	}
}
\theoremstyle{definition}

\begin{document}

\title{Evolutionary de Rham-Hodge method}

\author{Jiahui Chen$^{1}$,
Rundong Zhao$^{2}$,
Yiying Tong$^{2}$\footnote{
		Corresponding author,		Email: ytong@msu.edu},  ~  and  Guo-Wei Wei$^{1,3,4}$\footnote{
		Corresponding author,		Email: wei@math.msu.edu}\\
$^1$ Department of Mathematics,\\
Michigan State University, MI 48824, USA\\
$^2$ Department of Computer Science and Engineering,\\
Michigan State University, MI 48824, USA\\
$^3$ Department of Biochemistry and Molecular Biology,\\
Michigan State University, MI 48824, USA\\
$^4$ Department of Electrical and Computer Engineering\\
Michigan State University, MI 48824, USA
}

\date{}
\maketitle

\begin{abstract}
   The de Rham-Hodge theory is a landmark of the 20$^\text{th}$ Century's mathematics and has had a great impact on mathematics, physics, computer science, and engineering. This work introduces an evolutionary de Rham-Hodge method to provide a unified  paradigm for the multiscale geometric and topological analysis of evolving manifolds constructed from a filtration, which induces a family of evolutionary de Rham complexes. While the present method can be easily applied to close manifolds,  the emphasis is given to more challenging compact manifolds with 2-manifold boundaries, which require appropriate analysis and treatment of boundary conditions on differential forms to maintain proper topological properties. Three sets of unique evolutionary Hodge Laplacian operators are proposed to generate three sets of topology-preserving singular spectra, for which the multiplicities of zero eigenvalues correspond to exactly the persistent Betti numbers of dimensions 0, 1, and 2. Additionally, three sets of non-zero eigenvalues further reveal both topological persistence and geometric progression during the manifold evolution. Extensive numerical experiments are carried out via the discrete exterior calculus to demonstrate the utility and usefulness of the proposed method for data representation and shape analysis. 

\end{abstract}

\pagestyle{empty}

\newpage
{\setcounter{tocdepth}{3} \tableofcontents  }

\newpage
\pagestyle{plain}

\section{Introduction} \label{Introduction}
The de Rham-Hodge theory reveals that the cohomology of an oriented closed Riemannian manifold can be represented by harmonic forms.  It also holds for an oriented compact Riemannian manifold with boundary by forcing certain boundary conditions, such as absolute and relative cohomology~\cite{ray1971r}. This theory has been proved to be fundamentally important throughout algebraic geometry. It studies differential geometry and algebraic topology with partial differential equations (PDEs). The understanding of the de Rham-Hodge theory requires a variety of contemporary mathematical techniques including differential geometry, algebraic geometry, elliptic PDE, abstract algebra, topology, et al.

The de Rham-Hodge theory  has a wide range of applications, including not only mathematics, but also graphics/visualization~\cite{tong2003discrete, zhao20193d}, physics/fluids~\cite{foster1996realistic}, vision/robotics~\cite{gao2010singular,mochizuki2009spatial} and astrophysics/geophysics~\cite{mansour2004turbulence,akram2010regularisation}. Among all these applications, most of them rely upon the Hodge theory result, i.e., the Helmholtz-Hodge decomposition. It is one of the fundamental theorems in dynamical problems, describing a vector field into the gradient and curl components.

Due to the orthogonal decomposition, the analysis of vector fields becomes easier since certain properties such as incompressibility and vorticity of fluid dynamics can be studied on the orthogonal subspace. Such an orthogonal decomposition was first applied on a finite-dimensional compact manifold without boundary~\cite{hodge1989theory}, and then was developed for manifolds with boundaries~\cite{shonkwiler2009poincare}. Pushed by the visualization community, the implementation of orthogonal decomposition integrates a variety of boundary conditions with discrete vector fields expressed as discrete differential forms into two potential fields and harmonic fields~\cite{zhao20193d}.  The boundary conditions of the decomposition preserve orthogonality. The duality revealed by tangential and normal boundary conditions provides compact spectral representations of  the Laplace operators in the de Rham-Hodge theory. The spectra of de Rham-Laplace operators provide a quantitative approach to understanding topological spaces and geometry characteristics of manifolds and have been applied to biomolecular modeling and analysis \cite{zhao2019rham}. The development of discrete exterior calculus (DEC) is the driving force for de Rham-Hodge theory analysis and application \cite{arnold2006finite,desbrun2008discrete}. 

Over half a century ago, Kac asked a famous question, ``can one hear the shape of a drum?'' \cite{kac1966can}. Zelditch noticed that different drums may be distinguished by imposing restrictions with analytic boundary \cite{zelditch2000spectral}. However, the traditional spectral analysis cannot fully resolve the shape of a drum due to the isospectrum from different geometric shapes.  Innovative theoretical development is required to solve this long-standing spectral geometry problem. 

In the last few decades, geometric analysis has made great progress in understanding shapes that evolve in time. Geometric flows~\cite{willmore2013introduction} or geometric evolution equations have been extensively studied in mathematics~\cite{spruck1991motion,gomes2001using,mikula2004direct}, and many processes by which a curve or surface can evolve, such as the Gauss curvature flow and the mean curvature flow. Numerical techniques based on level sets were devised by Osher and Sethian~\cite{osher1988fronts} and have been extended and applied by many others in geometric flow analysis~\cite{wei2010differential,cecil2005numerical,du2004phase}. More recently, as the progress in contemporary life sciences, a large number of problems of unveiling the structure-function relationship of biomolecules and understanding of biomolecular systems, requires multiscale geometric modeling and analysis  ~\cite{bates2008minimal,wei2010differential,chen2011differential}.  However, compared with the investigations on curves and surfaces, a small amount of geometric explorations focuses on the evolution of compact manifolds specific to $\mathbb{R}^3$ due to the difficulty of computations. Additionally, it is rare to resolve topology from a nonlinear geometric PDE. Using a minimal molecular surface model \cite{bates2008minimal}, Wang and Wei studied the topological persistence via  the evolutionary profiles of the Laplace-Beltrami flow~\cite{wang2016object}. As a result, features of topological invariants are computed from the geometric PDE based filtration. In fact, there has been much effort in pure mathematics to understand the convergence of Riemannian manifolds in terms of sequences of submanifolds in metric spaces. However, the involved Gromov-Hausdorff distance can be computationally very  difficult. 

With the advancements in data  development and computational software, persistent homology has been promoted as a new multiscale approach for data analysis~\cite{zomorodian2005computing,edelsbrunner2010computational}. The traditional topological approaches describe the topology of a given object without invoking the metric or coordinate representations. Whereas, persistent homology bridges algebraic topology and multiscale analysis.  The essential difference is that persistent homology analyzes the persistence of the topological space through a filtration process, which is a family of simplicial complexes under a series of inclusion maps. Therefore a series of complexes is constructed based on filtration, which captures topological features changing over a range of spatial scales and reveals the features' topological persistence. In some sense, persistent homology can embed geometric information to topological invariants such that ``birth" and ``death" of connected components,  rings,  or cavities can be monitored by topological measurements during geometric scale changes. The original idea of varying scales was introduced by Frosini and Landi~\cite{frosini1999size} and by Robins in  1990s~\cite{robins1999towards}. Edelsbrunner et al. formulated the persistent homology and developed the first efficient computational algorithm \cite{edelsbrunner2000topological}. Zomorodian and Carlsson generalized the mathematical theory ~\cite{zomorodian2005computing}. Persistent homology has stimulated much theoretical development\cite{carlsson2009zigzag,edelsbrunner2010computational,chowdhury2018persistent,cang2018evolutionary,meng2019weighted, wang2019persistent}. Among them,  persistent spectral graph generates both topological persistence and spectral analysis \cite{wang2019persistent}.  Persistent homology has been applied to a variety of fields, including image analysis~\cite{carlsson2008local,pachauri2011topology,singh2008topological,bendich2010computing}, image retrieval~\cite{frosini2013persistent}, chaotic dynamics verification~\cite{mischaikow1999construction,kaczynski2006computational}, sensor network~\cite{de2005blind}, complex network~\cite{lee2012persistent,horak2009persistent}, data analysis~\cite{niyogi2011topological,wang2011branching}, computer vision~\cite{singh2008topological}, shape recognition~\cite{di2011mayer}, and computational biology~\cite{yao2009topological,xia2014persistent,xia2015persistent,gameiro2015topological,kovacev2016using}.

One of the first integrations of persistent homology and machine learning was developed for protein classification in 2015~\cite{cang2015topological}. Since then, persistent homology has been utilized as one of the most successful methods for the multiscale representation of complex biomolecular data~\cite{cang2017topologynet,cang2018integration,cang2018representability}. Two other multiscale representations of complex biomolecular data have also been proposed and found tremendous success in worldwide competitions in computer-aided drug design~\cite{nguyen2019mathematical,nguyen2019mathdl}. One of them is based on multiscale graphs~\cite{opron2015communication}, or more precisely, multiscale weighted colored graphs~\cite{bramer2018multiscale}. Eigenvalues of the graph Laplacians of multiscale weighted colored graphs were shown to provide some of the most powerful representations of protein-ligand binding interactions~\cite{nguyen2019agl}.   The other representation utilizes the curvatures computed from multiscale interactive molecular manifolds~\cite{nguyen2019dg}. The multiscale shape analysis offers an efficient means to discriminate similar geometries.  A common feature which is crucial to the success of the aforementioned three mathematical data representations is that they either create a family of multiscale topological spaces, or generate a family of multiscale graphs, or construct a family of manifolds, indicating the importance of the multiscale analysis in the representation of complex data with intricate internal structures.

Inspired by the aforementioned ideas, we introduce an evolutionary de Rham-Hodge method for data representation. The present evolutionary de Rham-Hodge method is developed by integrating differential geometry, algebraic topology, and multiscale analysis. It is noted that the fusion of algebraic topology and multiscale analysis leads to persistent homology, the combination of differential geometry and multiscale analysis renders manifold convergence \cite{sormani2010riemannian}, while the union of differential geometry and algebraic topology results in the de Rham-Hodge theory.  For a given dataset, using the evolutionary filtration developed in early work \cite{wang2016object}, we construct a sequence of evolving manifolds that lead to a geometry-embedded filtration under inclusion maps. The evolutionary de Rham-Hodge method is established on this sequence of manifolds. In general, the evolution of the manifolds can be either topological persistence which involves topological changes or geometric progression which does not involve topological changes.   We are interested in both the data analysis by evolutionary Hodge decompositions associated with various differential forms and the data representations via the evolutionary spectra of de Rham Laplace operators defined on the sequence of manifolds.  The evolutionary spectra reveal both the topological invariants and the geometric shapes of evolving manifolds. Such an evolutionary spectral analysis has great potential to ``hear the shape of a drum''.

In this work, we concern both close 2-manifolds and compact manifolds in $\mathbb{R}^3$  with boundaries, which require the enforcement of appropriate boundary conditions on differential forms to ensure topological properties.  Much effort has been given to the understanding and implementation of appropriate boundary conditions for the evolutionary de Rham-Hodge method, which results in three sets of unique evolutionary Hodge Laplacians. The multiplicities of the zero eigenvalues of these evolutionary  Hodge Laplacians provide the 0th, 1st, and 2nd persistent Betti numbers. Their non-zero eigenvalues further portray the geometric shape and topological characteristics of  data.

The rest of this paper is organized as follows. Section~\ref{primerDHT} is devoted to a brief review of the de Rham-Hodge theory, which includes the topics of the de Rham complex and Hodge decomposition. Then, the discrete forms and spectra generated by de Rham-Hodge theory are discussed in Section~\ref{primerDHT}. Readers familiar with the content in the above primer are recommended to start from Section~\ref{evoDHT}, where the evolutionary de Rham-Hodge method is formulated.  To demonstrate the utility and usefulness of the present method, we present the evolutionary de Rham-Hodge analysis of geometric shapes in Section~\ref{resDHT}. Finally, a conclusion is given in Section~\ref{conclusion}.

\section{A primer on de Rham-Hodge theory} \label{primerDHT}
To introduce the evolutionary de Rham-Hodge method, we briefly review the de Rham-Hodge theory to establish notation. We first discuss differential geometry and de Rham complex on smooth manifolds before reviewing the Hodge decomposition. Then, we illustrate the DEC discretization of the de Rham-Laplace operators and analyze their spectra.

\subsection{Differential geometry and de Rham complex}

Differential geometry is the study of shapes that can be represented by smooth manifolds of an arbitrary dimension.
A differential $k$-form $\omega^k \in \Omega^k(M)$ is an antisymmetric covariant tensor of rank $k$ on manifold $M$. Roughly speaking, at each point of $M$, it is a linear map from an array of $k$ vectors into a number, which switches sign if any two of the vectors are swapped. In general, it gives a uniform approach to define the integrals over curves, surfaces, volumes or higher-dimensional oriented submanifolds of $M$. More precisely, the antisymmetric rank-$k$ covariant tensor linearly maps $k$ edges from the first vertex of each $k$-simplex in a tessellation of the $k$-submanifold into a number, creating a Riemann sum that converges to an integral independent of the tessellation.

In $\mathbb{R}^3$, $0$-forms and $3$-forms can be recognized as scalar fields, as the antisymmetry permits one degree of freedom (DoF) per point, whereas $1$-forms and $2$-forms are considered vector fields as they require three DoFs per point. Our following discussion is specific to 3-dimensional (3D) volumes bounded by 2-manifolds in $\mathbb{R}^3$.

The \emph{differential} operator (i.e., exterior derivative) $d^k$ maps from the space of $k$-form on manifold, $\Omega^k(M)$ to $\Omega^{k+1}(M)$. It can be regarded as an antisymmetrization of the partial derivatives of a $k$-form. As such, it is a linear map $d^k:\Omega^k(M) \rightarrow \Omega^{k+1}(M)$ that satisfies the Stokes' theorem over any ($k\!+\!1$)-submanifold $\mathcal{S}$ in $M$:
\begin{equation}
\label{eqn:stokesThm}
\int_{\mathcal{S}} d^k \omega^k = \int_{\partial\mathcal{S}} \omega^k,
\end{equation}
where $\partial\mathcal{S}$ is the boundary of $\mathcal{S}$ and $\omega^k\in\Omega^k(M)$ is an arbitrary $k$-form.
Consequently, a key property of differential operator, $d^k d^{k-1}=0$, follows from that boundaries are boundaryless ($\partial\partial \mathcal{S}=0$). This implies that an \emph{exact} form (image of a $(k\!-\!1)$-form under differential) is \emph{closed} (i.e., is in the kernel of differential). The differential operator indeed provides a unification of a number of commonly used operators in 3D vector field analysis. Depending on the degree $k$ of differential forms, $d^k$ can be regarded as gradient ($\nabla$), curl ($\nabla\times$) and divergence ($\nabla\cdot$) operators for $0$-, $1$- and $2$-forms, respectively, e.g., $d^0$ takes the gradient of a scalar field (representing a $0$-form) to a vector field (representing a $1$-form).

With the linear spaces of $k$-forms treated as abelian groups under addition and the linear maps $d$ treated as group homomorphisms, they form a sequence that fits the definition of a \emph{cochain complex} as $d^k d^{k-1}=0$. This cochain complex of differential forms on a smooth manifold $M$ is known as the \emph{de Rham complex}:
\begin{center}
	\begin{tikzcd}
		0 \arrow{r} & \Omega^0(M) \arrow{r}{d^0} & \Omega^1(M) \arrow{r}{d^1} & \Omega^2(M) \arrow{r}{d^2} & \Omega^3(M) \arrow{r}{d^3} & 0.
	\end{tikzcd}
\end{center}
Note that $d^3$ maps $3$-forms to $4$-forms, but $k$-forms for $k>3$ are always zero in $\mathbb{R}^3$ due to antisymmetry.

The \emph{Hodge $k$-star} $\star^k$ (also called \emph{Hodge dual}) is linear map (and hence also a group isomorphism) from a $k$-from to its dual form, $\star^k : \Omega^k(M)\rightarrow\Omega^{n-k}(M)$. Due to the antisymmetry, both $k$-forms and their dual $(n\!-\!k)$-forms have the same DoF ${n \choose k}={n \choose n-k}$. More specifically, for an orthonormal basis $(\be_1,\be_2, \dots, \be_n)$, $\star^k(\be_{i_1}\wedge{\be_{i_2}}\wedge\cdots\wedge{\be_{i_k}})=\be_{j_1}\wedge{\be_{j_2}}\wedge\cdots\wedge{\be_{j_{n-k}}}$, where $\wedge$ denotes the antisymmetrized tensor product, and $(i_1,...,i_k,j_1,...,j_{n-k})$ is an even permutation of $\{1,2,...,n\}$. The associated $(\be_1,\be_2, \dots, \be_n)$ is a basis for $1$-forms, and $\be_{i_1}\wedge\dots\wedge\be_{i_k}$ form a basis for $k$-forms.

As $\star^k$ and $d^k$ can only operate on $k$-forms, we can omit the superscript of the forms or the operators when the dimension is clear from the context. The ($L_2$-)inner product of differential forms for two $k$-forms $\alpha, \beta \in \Omega^k(M)$ can be defined as
\begin{equation}
\label{eqn:innerProduct}
\langle \alpha , \beta \rangle = \int_M\alpha\wedge\star\beta = \int_M\beta\wedge\star\alpha.
\end{equation}
Under these inner products, the adjoint operators of $d$ are the \emph{codifferential} operators $\delta^k$: $\Omega^k(M)\rightarrow\Omega^{k-1}(M)$ , $\delta^k=(-1)^k\star^{4-k}d^{3-k}\star^k$ for $k=1,2,3$. In 3D, they can be identified with $-\nabla\cdot$, $\nabla\times$ and $-\nabla$ for $\delta^k$, $k= 1, 2, 3$ respectively in vector field analysis. Equipped with codifferential operators $\delta^k$, the spaces of differential forms now constitute a bi-directional chain complex,
\begin{center}
	\begin{tikzcd}
		\Omega^0(M) \arrow[r, shift left, "d^0"] &
		\Omega^1(M) \arrow[r, shift left, "d^1"] \arrow[l, shift left, "\delta^1"] &
		\Omega^2(M) \arrow[r, shift left, "d^2"] \arrow[l, shift left, "\delta^2"] &
		\Omega^3(M)                              \arrow[l, shift left, "\delta^3"].
	\end{tikzcd}
\end{center}
\begin{table}
	\small
	\caption{Exterior (odd rows) vs. traditional (even rows) calculus in $\mathbb{R}^3$. $f^0$, $\bv^1$, $\bv^2$ and $f^3$ stand for $0$-, $1$-, $2$- and $3$-forms with their components stored in either a scalar field $f$ or vector field $\bv$.}
	\begin{center}
		\renewcommand{\arraystretch}{1.3}
		\begin{tabular}{|c| l| l| l| l|}
			\hline 
			&	order $0$ & order $1$ & order $2$ & order $3$ \\
			\hline
			\hline
			form &  $f^0$ & $\bv^1(\ba)$ & $\bv^2(\ba,\bb)$ & $f^3(\ba,\bb,\bc)$ \\
			&  $f$ & $\bv\cdot \ba$ & $\bv\cdot(\ba\times\bb)$ &  $f[(\ba\times\bb)\cdot\bc]$ \\
			\hline
			$d$ & $df^0$ & $d \bv^1$ & $d \bv^2$ & $d f^3$ \\
			& $(\nabla f)^1$ & $(\nabla \times \bv)^2$ & $(\nabla \cdot \bv)^3$ &  $0$ \\
			\hline
			$\star$ &  $\star f^0$ & $ \star \bv^1$ & $\star \bv^2$ & $\star f^3$ \\
			& $f^3 $ & $\bv^2$ & $\bv^1$ & $f^0$ \\
			\hline
			$\delta$ & $\delta f^0$ & $\delta \bv^1$ & $\delta \bv^2$ & $\delta f^3$ \\
			& $0$ & $(- \nabla \cdot \bv)^0$ & $(\nabla \times \bv)^1$ &  $(-\nabla f)^2$ \\
			\hline
			$\wedge$ & $f^0\!\wedge\! g^0$ & $f^0\!\wedge\!\bv^1$ & $f^0\!\wedge\!\bv^2$, $\bv^1\!\wedge\!\bu^1$ & $f^0\!\wedge\!g^3$, $\bv^1\!\wedge\!\bu^2$\\
			& $(fg)^0$ & $(f\bv)^1$ & $(f\bv)^2$, $(\bv\!\times\!\bu)^2$ & $(fg)^3$, $(\bv\cdot\bu)^3$\\
			\hline
		\end{tabular}
	\end{center}
	\label{tb:ECoperator}
\end{table}
 Finally, the exterior calculus notations and their counterparts in traditional calculus are summarized in Table~\ref{tb:ECoperator}. The exterior calculus operations are strictly equivalent to the vector calculus operation in flat 3-dimensional space. A 0- or 3-form can be identified as a scalar function $f:M\subset\mathbb{R}^3\rightarrow\mathbb{R}$, while a 1- or 2-form is identified with a vector field $\bv:M\to \mathbb{R}^3$. Thus, we can use $f^0$, $\bv^1$, $\bv^2$ or $f^3$ to denote a scalar field $f$ or vector field $\bv$ regarded as a 0-, 1-, 2- or 3-form, respectively.

\subsection{Hodge decomposition for manifolds}
Hodge theory can be seen as the study of nonintegral parts (cohomology) of (scalar/vector) fields through the analysis of differential operators. Thus, it is often conveniently and concisely described by differential $k$-forms and the exterior calculus of these forms, as discussed in the previous section.

We first establish the aforementioned adjointness between the differential and codifferential operators.
Through integration by part and  the Stokes' theorem Eq.~(\ref{eqn:stokesThm}),
\begin{equation}
\label{eqn:integralPart}
\langle d\alpha, \beta\rangle=\langle \alpha, \delta\beta\rangle + \int_{\partial M} \alpha \wedge \star \beta.
\end{equation}
Thus, either for a boundaryless manifold ($\partial M=\emptyset$) or for forms that vanish on boundary ($\alpha|_{\partial M}=0$ or $\star\beta|_{\partial M} = 0$), the boundary integral vanishes, i.e., $\int_{\partial M} \alpha \wedge \star \beta = 0$. In such cases, the adjointness, $\langle d\alpha, \beta\rangle=\langle \alpha, \delta\beta\rangle$, implies that $d$ and $\delta$ satisfy the important property of adjoint operators---the kernel of a linear operator is the orthogonal complement of the range of its adjoint operator.

If we denote the space of \emph{normal forms} as $\Omega^k_n=\{\omega\in\Omega^k|\omega|_{\partial M} = 0\}$, and the space of \emph{tangential forms} as $\Omega^k_t=\{\omega\in\Omega^k|\star\omega|_{\partial M} = 0\}$, the orthogonal complementarity can be expressed as  $\Omega^k=\ker\delta^k \oplus d \Omega^{k-1}_n$ and $\Omega^k=\ker d^{k}\oplus \delta\Omega^{k+1}_t$. With $\im{d^{k-1}}\subset \ker d^k$ (based on the property of the cochain complex $d^k d^{k-1}=0$), the complementarity restricted to $\ker d^k$ implies
\begin{equation}
\label{eqn:kernelDecom}
\ker d^{k}=\mathcal{H}^k \oplus d\Omega^{k-1}_n,
\end{equation}
where $\mathcal{H}^k=\ker d^k\cap \ker \delta^k$ is the space of \emph{harmonic forms}, which are defined to be both closed and coclosed. Substituting the above equation into $\Omega^k=\ker d^{k}\oplus \delta\Omega^{k+1}_t$, we obtain the three-component Hodge decomposition,
\begin{equation}
\label{eqn:decomposition1}
\Omega^k = d\Omega_n^{k-1} \oplus \delta \Omega_t^{k+1} \oplus \mathcal{H}^k.
\end{equation}
Thus, any $\omega \in \Omega^k$ can be uniquely expressed as a sum of three $k$-forms
from the three orthogonal subspaces,
\begin{equation}
\label{eqn:orthogonalDiscription}
\omega = d\alpha_n + \delta\beta_t + h,
\end{equation}
where $\alpha_n \in \Omega_n^{k-1}$, $\beta_t \in \Omega_t^{k+1}$, and $h \in \mathcal{H}^k$. Note that the potentials $\alpha$ and $\beta$ do not have to be unique, and a variety of gauge conditions can be specified to make them unique.

\subsubsection{Boundaryless manifolds}

 When $\partial M\!=\!\emptyset$, $\Omega^k=\Omega^k_t=\Omega^k_n$, we can establish an isomorphism between the cohomology (of the de Rham complex described in the previous section) and the harmonic space, as was developed by Hodge.

In this case, Eq.~(\ref{eqn:kernelDecom}) can be written as
\begin{equation}
\label{eqn:kernelDecom1}
\ker d^{k}=\mathcal{H}^k \oplus \im d^{k-1}.
\end{equation}
Thus, we can find a unique element in $\mathcal{H}^k$ that corresponds to each equivalence class in the \emph{de Rham cohomology} $H^k_{dR}=\ker d^k/\im{d^{k-1}}$ (quotient spaces induced by the de Rham cochain complex). This bijection implies $\mathcal{H}^k \cong H^k_{dR}$, which indicates $\mathcal{H}^k$ is a finite-dimensional space with its dimension determined by the topology of the manifold.

Moreover, we can identify $\mathcal{H}^k$ as the kernel of a particular second-order differential operator, the \emph{de Rham-Laplace operator}, or \emph{Hodge Laplacian}, defined as $\Delta^k\equiv d^{k-1}\delta^k+\delta^{k+1} d^k$. Through the adjointness between $d$ and $\delta$, we have
\begin{equation}
\label{eqn:Laplace2form}
\langle \Delta \alpha, \alpha \rangle=\langle (d\delta+\delta d)\alpha, \alpha\rangle =  \langle d\alpha, d\alpha \rangle + \langle \delta\alpha, \delta\alpha \rangle.
\end{equation}
Denoting $\mathcal{H}^{k}_{\Delta} \equiv \ker\Delta^k$, the above equation implies that $\mathcal{H}^k_{\Delta} = \ker\Delta^k = \ker{d^k}\cap\ker{\delta^k}=\mathcal{H}^k$ for boundary-less manifolds.

As a direct consequence, we rewrite Eq.~(\ref{eqn:decomposition1}) as
\begin{equation}
\label{eqn:decomposition2}
\Omega^k = \im d^{k-1} \oplus \im \delta^{k+1} \oplus \mathcal{H}^k_\Delta.
\end{equation}
 The importance of the decomposition lies in that the first two components can be expressed as the derivatives of some potential functions, and the last non-integral part is spanned by the finite-dimensional harmonic space, whose dimension is determined by the topology of the domain due to the above-mentioned isomorphism. For example, for $\Omega^k$ with $k=1,2$, this decomposition is often recognized as the Helmholtz-Hodge decomposition of vector calculus in 3D, $\bv^1 = \nabla f^0 + \nabla \times \bu^2 + \bh^1$, and $\bv^2 = -\nabla f^3 + \nabla \times \bu^1 + \bh^2$.

\subsubsection{Manifolds with boundary}
For 3-manifolds with 2-manifold boundary, we need additional boundary conditions to have a finite dimensional kernel for the Laplacians, as in this case, $\mathcal{H}=\ker{d}\cap\ker{\delta}\subsetneq \mathcal{H}_{\Delta}$.
Through integration by part with the boundary, we have
\begin{equation}
\label{eqn:LaplacePart}
\langle \Delta \alpha, \alpha \rangle=\langle (d\delta+\delta d)\alpha, \alpha\rangle =  \langle d\alpha, d\alpha \rangle + \langle \delta\alpha, \delta\alpha \rangle+ \int_{\partial M} (\delta \alpha \wedge \star \alpha- \alpha \wedge \star d\alpha).
\end{equation}
Thus, if we can eliminate the boundary integral by restricting the space of forms, the kernel of $\Delta$ will be the intersection of the kernel of $d$ and $\delta$. Indeed, there are a variety of choices to satisfy boundary conditions, e.g., forcing the support of the differential form to be in the interior of manifolds. However, an option that is consistent with common physical boundary conditions is to restrict the differential form $\alpha$ in the decomposition
to be tangential to the boundary $\star\alpha|_{\partial M}=0$ or normal to the boundary $\alpha|_{\partial M}=0$ as we have required for the potentials.
%
Then, one natural choice to eliminate both terms in the boundary integral is to force $d \alpha$ to be tangential when $\alpha$ is tangential and force $\delta \alpha$ to be normal when $\alpha$ is normal. In other words, we modify the definition $\Omega_t$ to be the space of tangential forms with tangential differential, i.e., $\alpha_t\in\Omega_t$ if and only if
\begin{equation}
\label{eqn:tangentialBoundary}
\star \alpha_t |_{\partial M}=0, \quad \star d\alpha_t |_{\partial M}=0.
\end{equation}
Similarly, we modify the definition of $\Omega_n$ to be the space of normal forms with normal codifferential, i.e., $\alpha_n\in\Omega_n$ if and only if
\begin{equation}
\label{eqn:normalBoundary}
\alpha_n |_{\partial M}=0, \quad \delta\alpha_n |_{\partial M}=0.
\end{equation}


\begin{table}
	\caption{Boundary conditions of tangential and normal form}
	\begin{center}
		\renewcommand{\arraystretch}{1.3}
		\begin{tabular}{|l|l|l|l|l|}
			\hline
			type & $f^0$ & $\bv^1$ & $\bv^2$ & $f^3$ \\
			\hline \hline
			tangential & {\footnotesize unrestricted} & $\bv\cdot \bn = 0$ & $\bv \parallel \bn$ & $f|_{\partial M} = 0$ \\
			normal & $f|_{\partial M} = 0$ & $\bv \parallel \bn$ & $\bv \cdot \bn=0$ & {\footnotesize unrestricted} \\
			\hline
		\end{tabular}
	\end{center}
	\label{tb:boundary}
\end{table}

 To illustrate the boundary conditions explicitly, we consider a moving frame, which is formed at each boundary point by two tangent vectors of the boundary surface $\bt_1$ and $\bt_2$ and the normal vector to the surface $\bn$, with the typical convention that they form a right-hand orthonormal frame with the normal pointing outward. As a 1-form $\bv^1$ is tangential if $\star\bv^1(\bt_1,\bt_2)=\bv^2(\bt_1,\bt_2)= \bv\cdot(\bt_1\times \bt_2) = \bv\cdot\bn = 0$, it matches the condition that the corresponding vector field is tangential to the boundary. Similarly, a 1-form $\bv^1$ is normal to the boundary, if $\bv^1(\bt_i)=\bv\cdot\bt_i = 0$ for $i=1,2$, thus it is the equivalent to $\bv$ is normal to the boundary. For a $2$-form $\bv^2$, its normal (tangential) boundary condition is the same as the tangential (normal) boundary condition of $\bv^1$. Therefore, normal (tangential) $2$-forms should have their corresponding vector fields tangential (normal, resp.) to the boundary. Additionally, tangential $3$-forms (normal $0$-forms) are zero on the boundary whereas normal $3$-forms (tangential $0$-forms) automatically satisfy the boundary condition. In Table~\ref{tb:boundary}, we summarized these choices of the boundary conditions for tangential and normal $k$-forms in 3D.

In vector field representation, the boundary conditions Eqs.~(\ref{eqn:tangentialBoundary}) and (\ref{eqn:normalBoundary}) are equivalent to the following. The choice of a $1$-form in $\Omega^1_t$ (a $2$-form in $\Omega^2_n$) is equivalent to enforcing a tangential vector field $\bv$ to have its curl to be normal to the boundary, i.e., adding two homogeneous Neumann boundary conditions to the (Dirichlet-type) tangentiality,
\begin{equation}
\label{eqn:bdyT1form}
\bv\cdot\bn=0, \quad \nabla_{\bn}(\bv\cdot\bt_1)=0,\quad \nabla_{\bn}(\bv\cdot\bt_2)=0.
\end{equation}
For a normal vector field $\bv$ ($1$-forms in $\Omega^1_n$ or $2$-forms in $\Omega^2_t$), it amounts to adding one homogeneous Neumann boundary condition derived from the zero divergence on the boundary to the  (Dirichlet-type) orthogonality constraints,
\begin{equation}
\label{eqn:bdyT2form}
\bv\cdot\bt_1=0,\quad \bv\cdot\bt_2=0,\quad \nabla_{\bn}(\bv\cdot\bn)=0.
\end{equation}
For an unrestricted function $f$ (tangential $0$-forms or normal $3$-forms), it amounts to forcing its gradient to be tangential at the boundary (Neumann-type),
\begin{equation}
\label{eqn:bdyT0form}
\nabla_{\bn} f|_{\partial M}=0,
\end{equation}
and a function $f$ for tangential $3$-forms (normal $0$-forms) satisfies the homogeneous Dirichlet boundary condition
\begin{equation}
\label{eqn:bdyT3form}
f|_{\partial M}=0.
\end{equation}
 
With these modified boundary conditions, we still have the same Hodge decomposition,
\begin{equation}
\label{eqn:decomposition4}
\Omega^k = d\Omega_n^{k-1} \oplus \delta \Omega_t^{k+1} \oplus \mathcal{H}^k.
\end{equation}
This is due to the fact that $d\Omega_n$ (or $\delta\Omega_t$) remains the same regardless of whether $\Omega_n$ (or $\Omega_t$) contains the additional boundary conditions, as they can be seen as part of the gauge condition that restricts the potentials but not their differential (codifferential).

As mentioned above, with the boundary, $\mathcal{H}^k$ is no longer finite dimensional or the kernel of of Laplacians $\mathcal{H}^k_\Delta$. However, if we restrict $\Delta$ to $\Omega_t$ or $\Omega_n$ and denote the corresponding operator as $\Delta_t$ and $\Delta_n$ respectively, we can still find finite dimensional kernels $\mathcal{H}^k_{\Delta_t}$ and $\mathcal{H}^k_{\Delta_n}$ that correspond to $\mathcal{H}^k\cap \Omega_t$ or $\mathcal{H}^k\cap \Omega_n$ orthogonal to $\im d$ and $\im \delta$.

 In fact, the harmonic space $\mathcal{H}^k$ can be further decomposed into tangential, normal harmonic forms and exact-coexact harmonic forms $\mathcal{H}^k=(\mathcal{H}^k_{\Delta_t} + \mathcal{H}^k_{\Delta_n}) \oplus (d\Omega^{k-1}\cap \delta\Omega^{k+1})$ as proposed by Friedrichs~\cite{friedrichs1955differential}.  Moreover, in flat 3D space, all three subspaces are orthogonal to each other. The third space can be seen as the infinite-dimensional space of solutions to Laplace equations in dimension $k\pm 1$ with either normal or orthogonal boundary conditions. Thus, we can focus on the Laplacian operators that are either tangential or normal for analysis.

In total, there are 8 different Hodge Laplacians ($\Delta^k_{t}$ and $\Delta^k_{n}$ for $k=0,1,2,3$)
and 8 associated finite dimensional harmonic spaces. Friedrichs also noted that for manifolds with boundary, the tangential harmonic spaces are isomorphic to the absolute de Rham cohomology $\mathcal{H}^k_{\Delta_t} \cong H^k(M)$, and the normal harmonic spaces are isomorphic to the relative de Rham cohomology $\mathcal{H}^k_{\Delta_n} \cong H^k(M,\partial M)$. From the dimensionality of the corresponding homology (Betti numbers) of the manifold $M$, together with the Hodge duality between $\mathcal{H}^k_{\Delta_t}$ and $\mathcal{H}^{3-k}_{\Delta_n}$,
we can obtain the dimensions of all these harmonic spaces: $\beta_k=\dim \mathcal{H}^k_{\Delta_t} = \dim \mathcal{H}^{3-k}_{\Delta_n}$. Roughly, speaking, $\beta_0$ is the number of connected components, $\beta_1$ is the number of rings, $\beta_2$ is the number of cavities, and $\beta_3$ is 0 as $M$ in flat 3D cannot contain any noncontractible topological 3-sphere.

\subsection{Discrete forms and spectral analysis} \label{Discreteformandspectrumanalysis}

 In practical applications, the de Rham-Hodge theory is often computed for decompositions and spectral analysis. In both cases, the discretization of exterior derivatives is required.   We follow one typical discretization of the exterior calculus on differential forms, the discrete exterior calculus (DEC)~\cite{desbrun2008discrete}. A major technical aspect is the handling of arbitrarily complex geometric shapes in 3D. In spectral analysis, the Hodge Laplacian operators and their boundary conditions are to be implemented such that the key topological property of   $d\circ d=0$, which defines the de Rham cohomology,  is preserved in the discrete version by DEC in complex computational domains. First, the domain of differential forms, in this case, a 3-manifold embedded in 3D Euclidean space is tessellated into a 3D simplicial complex, i.e.,  a tetrahedral mesh. Any $k$-form $\omega$ is represented by its integral on oriented $k$-D elements ($k$-simplex) of the mesh, listed as a vector $W$ with the length equaling the number of $k$-simplices. More specifically, a discrete $0$-form is the assignment of one real number per vertex, a discrete $1$-form is the assignment of one value per oriented edge, a discrete $2$-form is the assignment of one value per oriented triangle, and a discrete $3$-form is the assignment of one value per tetrahedron (tet). The choice of orientation per $k$-simplex is arbitrary since the antisymmetry of a $k$-form guarantees that the integral on that $k$-simplex only changes its sign. 

Now the linear operator $d^k$ is represented by a sparse matrix ${D}_k$, which is implemented as the transpose of the signed incidence matrix between $k$-simplices and $(k\!+\!1)$-simplices, with the sign determined by mutual orientation. This can be seen as the consequence of the aforementioned Stokes' theorem, because the integral of $d\omega$ on each $(k\!+\!1)$-simplex is exactly the sum of the integral of $\omega$ on the boundary of the $(k\!+\!1)$-simplex, which is the union of its consistently oriented $k$-simplex faces.

Thus, the defining property in de Rham-Hodge theory ${D}_{k+1}{D}_{k}=0$ is preserved through as the boundary of the boundary is empty. As shown in Fig.~\ref{fig:discretedR}, the adjoint operator $\delta^k$ is implemented as ${S}_{k-1}^{-1}{D}^T_{k-1}{S}_{k}$, where ${S}_k$ is discretization of the $L_2$-inner product between two discrete $k$-forms such that $(W_1^k)^TS_k W_2^k$ is an approximation of $\langle\omega_1^k,\omega_2^k\rangle$. In this work, we use the lowest order diagonal matrices for $S_k$ for simplicity, but higher-order Galerkin matrices for  $k$-form basis can be developed with proper treatment on matrix inversion for better accuracy. Such a discrete Hodge star operator can also be seen as a mapping from a discrete $k$-form to a discrete dual $(3\!-\!k)$-form defined on the basis associated with dual elements of a dual mesh to the tet mesh. Obviously, this field needs more effort from the computational mathematics  community  

With both the differential operators and the Hodge stars discretized, the discrete counterpart of a Hodge Laplacian $\Delta^k$ is defined as $S_k^{-1}L_k$ through products and summations of these matrices following the continuous version, here 
\begin{equation}
\label{eqn:discreteLaplacian}
L_k=D_k^TS_{k+1}D_k+S_kD_{k-1}S^{-1}_{k-1}D^T_{k-1}S_k.
\end{equation}
The reason that  $L_k$ is used frequently as the discrete Hodge Laplacian instead of $S_k^{-1}L_k$ is its symmetry. Alternatively, we can also see $L_k$ as the quadratic form on the space of discrete $k$-forms, such that $W^T L_k W$ is an approximation of $\langle \omega, \Delta \omega\rangle$.

In our analysis of volumetric shapes, we conjecture that the evolution of topological and geometric structures is related not only to the null spaces of Hodge  Laplacians, but also to the general spectra of these operators, in particular, those eigenvalues that are close to zero. The associated eigen differential forms can be found through a generalized eigenvalue problem for the discrete Hodge Laplacian and Hodge star operators.
\begin{equation}
\label{eqn:genEigval}
L_k W^k=\lambda^k S_k W^k.
\end{equation}

\begin{figure}
	\begin{center}
		\begin{tikzcd}
			\text{0-form} \arrow[r, "{D}_0"] \arrow[d, xshift=0.7ex, "{S}_0"] &
			\text{1-form} \arrow[r, "{D}_1"] \arrow[d, xshift=0.7ex, "{S}_1"] &
			\text{2-form} \arrow[r, "{D}_2"] \arrow[d, xshift=0.7ex, "{S}_2"] &
			\text{3-form}                         \arrow[d, xshift=0.7ex, "{S}_3"] \\
			\text{dual 3-form} \arrow[u, xshift=-0.7ex, "{S}_0^{-1}"] &
			\text{dual 2-form} \arrow[u, xshift=-0.7ex, "{S}_1^{-1}"] \arrow[l, "D_0^{T}"] &
			\text{dual 1-form} \arrow[u, xshift=-0.7ex, "{S}_2^{-1}"] \arrow[l, "D_1^{T}"] &
			\text{dual 0-form} \arrow[u, xshift=-0.7ex, "{S}_3^{-1}"] \arrow[l, "D_2^{T}"]
		\end{tikzcd}
		\caption{Discrete de Rham cohomology; $D_k$ is the combinatorial operators such that
			${D}_{k+1}{D}_{k}=0$; ${S}_k$ is the discrete Hodge stars.}
		\label{fig:discretedR}
	\end{center}
\end{figure}
For illustration purpose, we can reformulate Eq.~(\ref{eqn:genEigval}) as a regular eigenvalue problem,
\begin{equation}
\label{eqn:genRegular}
\bar{L}_k\bar{W}^k=\lambda^k \bar{W}^k,
\end{equation}
where $\bar{L}_k = S_k^{-1/2}L_kS_k^{-1/2}$ and $\bar{W}^k=S_k^{1/2} W^k$.
Then, to partition the spectrum of the modified discrete Hodge Laplacian, we express it as the sum of two semi-positive-definite matrices,
\begin{equation}
\label{eqn:modifiedLaplacian}
\bar{L}_k=\bar{D}_k^T\bar{D}_k+\bar{D}_{k-1}\bar{D}_{k-1}^T,
\end{equation}
where $\bar{D}_k=S^{1/2}_{k+1}D_kS_k^{-1/2}$. We can observe that the cohomology structure is maintained as $\bar{D}_{k+1}\bar{D}_k=0$. Moreover, now the adjoint operator of $\bar{D}_k$, in the $L_2$ inner products defined by the Hodge stars, is simply its transpose $\bar{D}^T_k$. Thus, the entire spectrum of $\bar{L}_k$ can be studied through the singular value decomposition of the discrete differential operator
\begin{equation}
\label{eqn:SVDdisdiffoperator}
\bar{D}_k=U_{k+1}\Sigma_kV_k^T,
\end{equation}
where $U_{k+1}$ and $V_k$ are orthogonal matrices, and $\Sigma_k$ is a rectangular diagonal matrix with non-negative real elements. We can recognize the nonzero spectra of the modified Hodge Laplacian as the union of the squares of the nonzero entries from $\Sigma_k$ and $\Sigma_{k-1}$, since
\begin{equation}
\label{eqn:SVDmodifiedLaplacian}
\bar{L}_k=V_k\Sigma_{k}^2V^T_k+U_k\Sigma_{k-1}^2U^T_k.
\end{equation}
Note that for $0$- or $3$-forms, one of the $\Sigma$'s contains only zeros.

Based on the Hodge decomposition Eq.~(\ref{eqn:decomposition4}), we can also notice that the columns of $V_k$ that correspond to nonzero singular values in Eq.~(\ref{eqn:SVDmodifiedLaplacian}) are orthogonal to those of $U_k$, which means the entire $k$-form space is spanned by harmonic forms (eigen form with eigenvalue 0), and those column vectors of $V_k$ and $U_k$.

For domains with boundaries, the tangential or normal forms are restricted by Dirichlet and/or Neumann boundary conditions, which can be implemented by whether to include the boundary elements or not for $D_k$.
We denote the discrete differential operator for tangential (normal) $k$-forms as $D_{k,t}$ (respectively $D_{k,n}$).
For the detail on the construction of these matrices, readers are referred to our previous work~\cite{zhao20193d}.
In summary, for the four types of $k$-form ($k=0,1,2,3$) with two boundary conditions, there are 8 different discrete Hodge Laplacians ($L_{k,t}$ and $L_{k,n}$) in total, such that
\begin{equation}
\label{eqn:Laplacianktn}
\begin{split}
L_{k,t}&=D_{k,t}^TS_{k+1}D_{k,t}+S_kD_{t,k-1}S^{-1}_{k-1}D^T_{t,k-1}S_k,\\
L_{k,n}&=D_{k,n}^TS_{k+1}D_{k,n}+S_kD_{n,k-1}S^{-1}_{k-1}D^T_{n,k-1}S_k.
\end{split}
\end{equation}
Based on the above singular value analysis, the non-zero spectrum of $\bar{L}_k$ is the union of squared singular values of $\bar{D}_k$ and those of $\bar{D}_{k-1}$. Therefore, for each type of boundary conditions, the spectra of the four discrete Hodge Laplacians only depend on the singular spectra of $\bar{D}_0$, $\bar{D}_1$ and $\bar{D}_2$. Furthermore, in Table~\ref{tb:boundary}, the same set of boundary conditions is shared between tangential $1$-forms and the normal $2$-forms, between tangential $2$-forms and normal $1$-forms, between normal $3$-forms and tangential $0$-forms, and between tangential $3$-forms and normal 0-forms. This duality between tangential $k$-forms and normal ($3\!-\!k$)-forms is also present in the corresponding operators between these forms, more specifically, the equivalence exists between $\bar{D}_{0,t}$ and $\bar{D}_{2,n}^T$,
$\bar{D}_{1,t}$ and $\bar{D}_{1,n}^T$, and
$\bar{D}_{2,t}$ and $\bar{D}_{0,n}^T$.
We thus reduce the 8 different spectra of Hodge Laplacians to 3 distinct sets of different singular spectra. We denote the set of singular values of $\bar{D}_{0,t}$ for the tangential gradient eigen field by $T$, the set of the singular values of $\bar{D}_{1,t}$ for the curl eigen field by $C$, and the set of the singular value set of $\bar{D}_{2,t}$ for tangential divergent eigen field by $N$.

Although each of the 8 spectra for Hodge Laplacians defined on smooth manifolds can be represented by the combination of one or two sets of the $T$, $C$ and $N$, the numerical calculations of the singular values of the equivalent differential operators can deviate from these due to the different DoFs in the representations for different discrete forms, as well as the inaccuracy introduced by the approximation of Hodge star and differential operators. While the numerically computed singular values of tangential $k$-forms $\bar{D}_{k,t}$ can deviate from those of normal $(3\!-\!k)$-forms $\bar{D}_{2-k,n}^T$, as the observation in previous work~\cite{zhao2019rham}, with increased resolution, the low frequencies converge reasonably well.

\section{Evolutionary de Rham-Hodge method}\label{evoDHT}

  In this section, we introduce the evolutionary de Rham-Hodge method to analyze the topological and geometric properties throughout the evolution of manifolds.  We first discuss the existing data that motivates the present theoretic formulation. Then, we provide the mathematical description of manifold evolution, followed by the definitions of the associated persistence and progression. We extend the usual study of cohomology (associated to zero eigenvalues of Hodge Laplacians) to employing the leading small non-zero eigenvalues to facilitate the concepts of persistence and progression so that the variations of topological spaces ($\beta_0$, $\beta_1$ and $\beta_2$) can be traced to the changes in the eigenvalues away from or towards zero as the geometry evolves.

\subsection{Data and their de Rham-Hodge analysis}
 
 Most commonly occurred data are closed manifolds, such as star surfaces,  earth surfaces, brain surfaces, and molecular surfaces. The de Rham Laplace operator can be applied to compute eigenfunctions and eigenvalues for the geometric shape analysis. Another interesting type of data includes scalar or vector functions defined on closed manifolds, such as temperature or ocean currents on the earth's surface and in compact manifolds with boundaries, such as the electron densities or electrostatic potentials in proteins or the magnetic fields around the earth.  The Hodge decomposition can be directly applied to these functions.  For smooth scalar functions, surface contours can be specified to generate compact manifolds with boundaries. The geometric shape analysis via the de Rham Laplace operator can be carried out. A special class of data is the density distributions, either obtained from cryogenic electron microscopy (cryo-EM),   magnetic resonance imaging (MRI) or created from quantum mechanical calculations. In this situation, one can render a family of  inclusion surfaces by systematically varying the density isovalues. The de Rham-Hodge analysis and modeling of this family of  inclusion surfaces are the objects of the present theoretical development.  

 The evolutionary de Rham-Hodge method developed in this work can also be applied to point cloud data, such as stars in the universe, atoms in biomolecules,  and the output of 3D scanning processes. In this situation, one can carry out a discrete to continuum map to create volumetric density functions from point clouds \cite{xia2014multiscale,nguyen2016generalized}. Then, a  family of inclusion surfaces can be obtained for the evolutionary de Rham-Hodge analysis. 
     
Flexibility rigidity index (FRI) density is a useful tool to construct a continuous density distribution from a set of discrete point cloud data inputs.  By selecting an isovalue  from the FRI density, one can further generate a boundary surface, which composes the 3-manifold with a 2-manifold boundary. FRI density has been shown to be particularly straightforward to implement and computationally stable on any point cloud~\cite{nguyen2016generalized} and is defined by the following position-dependent rigidity (or density) function~\cite{xia2014multiscale}
\begin{equation}
\label{eqn:densityFcn}
\rho(\br,  \eta) = \sum_{j=1}^{N}\Phi(\|\br-\br_j\|; \eta)
\end{equation}
where $\br$ is a point in space, $N$ is the number of particles, $\br_j$ is the location of a data point $j$, $\eta$ is a  scaling parameter, and $\Phi(\cdot;\eta)$ is a correlation function, i.e., a real-valued monotonically decreasing function with the following admissibility conditions
\begin{equation}
\label{eqn:densityCond}
\begin{split}
\Phi(\|\br-\br_j\|;\eta) &= 1, \text{ as } \|\br-\br_j\| \rightarrow 0,\\
\Phi(\|\br-\br_j\|;\eta) &= 0, \text{ as } \|\br-\br_j\| \rightarrow \infty,
\end{split}
\end{equation}
One used families of correlation functions is the generalized exponential functions
\begin{equation}
\label{eqn:expFcn}
\Phi(\|\br-\br_j\|;\eta) = \exp({-(\|\br-\br_j\|/\eta )^\kappa}), \quad \kappa>0.
\end{equation}
Here, the weight $\eta$ is application-dependent, e.g., the multiplication of a scaling parameter and the van der Waals radius $r_{{\rm vdw}_j}$ of the atom at $\br_j$ for molecular data. In fact, $\eta$ can be chosen as anisotropic function to induce a multidimensional persistent homology filtration \cite{xia2015multidimensional}. In our numerical tests, we use the generalized exponential function with $\kappa=2$, which is known as the Gaussian function. A family of 3-manifolds can be defined by a varying level set parameter (isovalue) $c\in(0,c_\text{max})$, where $c_\text{max}=\max{\rho(\cdot, \eta)}$,
\begin{equation}
\label{eqn:manifolds}
M_c=\{\br|\rho(\br, \eta) \leq c_\text{max}-c\},
\end{equation}
which has the level-set of $\rho$ as its boundary $\partial M_c = \{\br|\rho(\br, \eta) = c_\text{max}-c\}$.

\subsection{Manifold evolution}

Hodge theory studies the de Rham cohomology groups of a smooth manifold $M$, and established the bijection from equivalence classes in a cohomology group to a harmonic differential form in the null space of the corresponding Hodge   Laplacian. While these harmonic forms associated with the zero eigenvalues in the spectra of Hodge Laplacians carry some geometric information in addition to the topology, the non-zero spectra provide richer geometric information than the multiplicity of zero. However, the geometry is not uniquely determined by the spectra of the Hodge Laplacians (even for planar shapes), as one cannot hear the shape of a drum~\cite{kac1966can}. Thus, we propose to extend the study of de Rham-Hodge theory to a family of smooth manifolds instead of one specific manifold and track the spectral changes in a sequence of manifolds. Such a family of manifolds controlled by a continuous filtration parameter is sometimes called the evolution of manifolds embedded in an ambient manifold, which in our case is the 3D Euclidean space.
 
 The evolution of manifolds is often defined through a smooth map from a basic manifold $B$ to a family of submanifold $\{M_c\}$ of an ambient manifold $M$ at a given instant (the value of parameter $c$ treated as time). More precisely, it is the smooth map $F: B \times [0,c_\text{max}] \rightarrow M$ such that $F^c=F(\cdot, c)$ is an immersion for every $c$. The one-parameter family of subsets of $M$, $\{F^c(B)\}_{c\geq 0}$  is then called the \emph{evolving manifold}. However, such a Hodge Lagrangian description makes it hard to handle topological changes, especially if each mapping is restricted to be an embedding. Therefore, in this work, we directly use the Eulerian representation described by $M_c$ in Eq.~(\ref{eqn:manifolds}). This level-set bounded volume evolution handles both the geometric progression and topological changes in a consistent fashion. As Morse functions are dense in continuous functions, we can assume $\rho(\br,\eta)$ to be a Morse function without loss of generality, since otherwise, we can use symbolic perturbation to make it a Morse function. We can regularly sample the interval $(0,c_\text{max})$ at $n$ sample locations, forming an index set $I=\{c_0,c_1,...,c_n\}$, such that none of the parameters are one of the isolated critical values through symbolic perturbation if necessary. Noting that $M_c$ are only non-manifold when $c$ is a critical point of the Morse function, the snapshots of the evolving manifold, $\{F^c\}_{c\in I}$, are all manifolds. Thus, they form a filtration of manifold $M$, with the inclusion map $\mathfrak{I}_{l,l+1}: M_{l}\hookrightarrow M_{l+1}$ linking each pair of consecutive manifolds and

\begin{center}
	\begin{tikzcd}
		M_{0} \arrow[r, "\mathfrak{I}_{0,1}"]  &
		M_{1} \arrow[r, "\mathfrak{I}_{1,2}"] &
		M_{2} \arrow[r, "\mathfrak{I}_{2,3}"] &
		\cdots \arrow[r, "\mathfrak{I}_{n-1,n}"] &
		M_{n} \arrow[r, "\mathfrak{I}_{n,n+1}"] &
		M=M_{c_\text{max}}.
	\end{tikzcd}
\end{center}

If $(c_l,c_{l+p})$ does not contain any critical points of $\rho(\br, \eta)$ and the largest critical value smaller than $c_l$ is $c_c$, the inclusion map $\mathfrak{I}_{l,l+p}:M_l\hookrightarrow M_{l+p}$ is also homotopic to a homeomorphism from $M_{l}$ to $M_{{l+p}}$, which can be constructed by moving every point $\br$ with $\rho(\br,\eta)>c_{\max}-c_c$ along the gradient integral line of $\rho(\cdot, \eta)$ to a point $\hat\br$ such that {$\rho(\br,\eta)-\rho({\hat\br}, \eta)=(c_{l+p}-c_l)e^{1-\frac{c_l-c_c}{\rho(\br, \eta)-c_{\max}-c_c}}$}. When the two parameter values are similar, one can also see that the above map is nearly isometric since the deformation is close to an identity map.

When $(c_l,c_{l+p})$ contains critical points of the Morse function, there is no smooth homeomorphism between $M_l$ and $M_{l+p}$ as the level set underwent topological changes. Without loss of generality, we can assume that there is only one critical point, which can be classified as (local) minimum, 1-saddle, 2-saddle, or (local) maximum, based on the signature of the Hessian of $\rho$. As all minima of $\rho$ is at the value of $0$, the interval may only contain the latter three types: if it is a maximum, one 2nd homology generator in $M_l$ will be mapped to $0$ in $M_{l+p}$ for the mapping induced by the inclusion; if it is a 2-saddle, either $M_l$ has a 1st homology generator mapped to $0$ or $M_{l+p}$ contains a 2nd homology generator not in the image of the induced mapping from $H(M_l)$ to $H(M_{l+p})$; similarly, if it is a 1-saddle, either $M_l$ has a 0th homology generator mapped to $0$ or $M_{l+p}$ contains a 1st homology generator, not in the image of the induced mapping. Through the isomorphisms among the de Rham cohomology, singularly homology, simplicial homology, and simplicial cohomology, we can use the persistent homology to study the mapping between the de Rham cohomologies indirectly. However, we found that direct construction can reveal some additional insight on the relation and persistence of the harmonic forms across different manifolds, as we discuss next.

\subsection{Persistence of harmonic forms}

\subsubsection{Normal harmonic forms}

Drawing an analogy from persistent homology, we first attempt to construct a homomorphism from closed forms on $M_l$ to closed forms on $M_{l+p}$, i.e., from $\ker d_{l}$ to $\ker d_{l+p}$, if we use the subscript $l$ to denote the operator defined on $M_l$. For manifolds with boundary, one realizes that this is not possible for tangential forms through the isomorphism relations to cochain and chain spaces on simplicial complexes, but rather straightforward for normal forms in the discrete case. More specifically, we can map $k$-forms in $M_l$ by setting values for simplices in $M^c_{l,p}=M_{l+p}\backslash M_l$ to $0$, i.e., a 0-padded $k$-cochain on $M_{l+p}$ as the image of a $k$-cochain on $M_l$ assuming that $M_l$ has a tessellation that is a subcomplex of the tessellation of $M_{l+p}$. The reason that the image of $\omega_l \in \ker d_{l}$ remains in $\ker d_{l+p}$ is that the value of $d\omega_{l+p}$ on any $({k\!+\!1})$-simplex with one or more faces in $\partial M_l$ is still $0$, as $\omega_{l}|_{\partial M_l} = 0$.

However, in the continuous case, setting $\omega$ to $0$ in $M^c_{l,p}$ creates either discontinuity or at least large $\delta\omega$ near the boundary. A smoother extension of the $\omega$ from $M_l$ to $M_{l+p}$ can be defined by minimizing the Dirichlet energy $\langle d\omega,d\omega\rangle + \langle \delta \omega,\delta \omega\rangle$ in $M^c_{l,p}$, which leads to simply a Laplace equation $\Delta \omega = 0$. The boundary of $M^c_{l,p}$ is the union of $\partial M_{l}$ and $\partial M_{l+p}$ with the orientation of the former flipped. Recall that when $\omega$ is normal to the boundary i.e., $\omega_l|_{\partial M_l}=0$, we also impose the condition that $\delta \omega$ is normal to the boundary ($\delta\omega_l|_{\partial M_l}=0$). For the extension, we keep this condition on $\partial M_{l+p}$, while on $\partial M_l$ we impose the continuity instead, $\omega_{l+p}|_{\partial M_l} = \omega_l|_{\partial M_l}$. Note that the resulting Laplace equation has a finite kernel identical to that of $\Delta_n$ on $M^c_{l,p}$, so we can find a unique solution by forcing the solution to have 0 projection to this kernel~\cite{zhao20193d}.

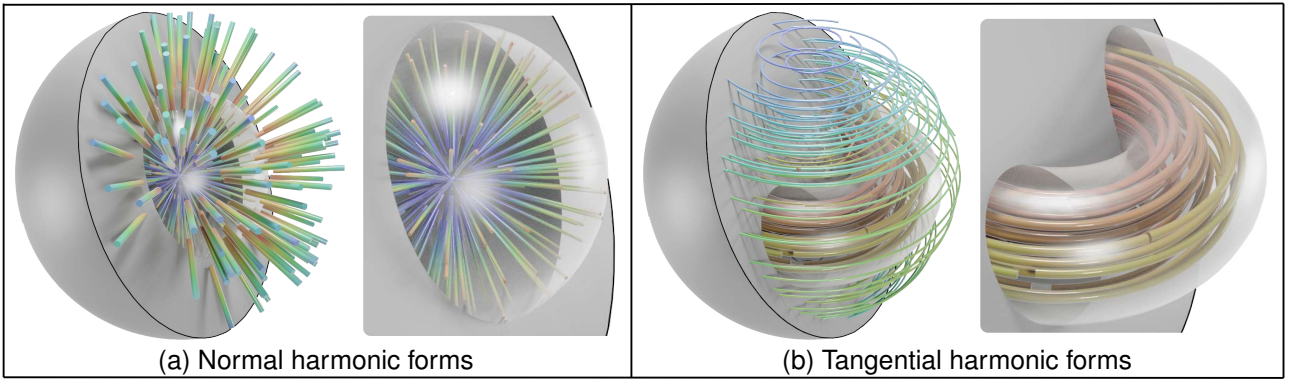
\begin{figure*}[t!]
	\centering
	\includegraphics[height=2in]{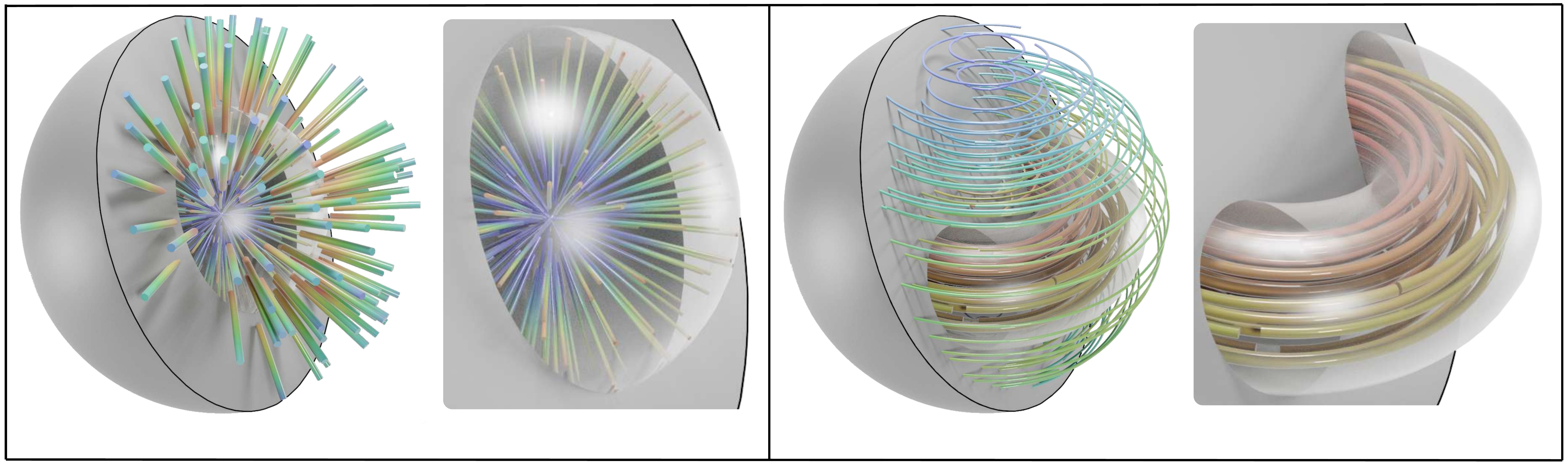}
	\begin{picture}(17,0)
	\put( -175, 17){\text{(a) Normal harmonic forms}}
	\put(   61, 17){\text{(b) Tangential harmonic forms}}
	\end{picture}
	\caption{Illustration of normal and tangential harmonic field extensions. Thick lines are the inputs and thin lines are the extended outputs. Left charts in both (a) and (b) show harmonic fields and their extensions while right charts give meticulous detail of interior parts. (a) Normal harmonic forms. A solid ball with a cavity extends inward to a solid ball without cavity. The outside surface is fixed. (b) Tangential harmonic forms. A torus extends to a solid ball.}
	\label{fgr:normalTangentialVF}
\end{figure*}

For instance, if we have a normal 1-form $\omega_l$ to extend, we can impose the homogeneous boundary condition for the proxy vector field $\bv$ on $\partial M_{l+p}$ as in Eq.~(\ref{eqn:bdyT2form}),
\begin{equation}
\label{eqn:bc_ext_n}
\bv_{l+p}\cdot\bt_1=0,\quad \bv_{l+p}\cdot\bt_2=0,\quad \nabla_{\bn}(\bv_{l+p}\cdot\bn)=0;
\end{equation}
whereas on $\partial M_l$, we use a Dirichlet boundary condition for continuity $\bv_{l+p}=\bv_l$, i.e.,
\begin{equation}
\label{eqn:bc_ext_n_1}
\bv_{l+p}\cdot\bn=\bv_{l}\cdot\bn, \quad \bv_{l+p}\cdot\bt_1=0,\quad \bv_{l+p}\cdot\bt_2=0.
\end{equation}
We denote the map through this harmonic extension as $\mathfrak{E}_{l,p}$, i.e., $\omega_{l+p} = \mathfrak{E}_{l,p}(\omega_l)$. However, the minimization of Dirichlet energy does not imply $\delta \omega_{l+p} = 0$ even when $\delta \omega_l = 0$. Nevertheless, $d \omega_{l+p} = 0$ is always possible, since otherwise, one would be able to perform a Hodge decomposition to find a tangential $(k\!+\!1)$-form $\beta_t$ in $M^c_{l,p}$ and remove $d\omega_{l+p}$ by subtracting $\delta\beta_t$ from $\omega_{l+p}$. An alternative is to restrict the extension to minimize $\langle \delta\omega,\delta\omega\rangle$ under the constraint $d\omega_{l+p}=0$ in $M^c_{l,p}$, which results in a fourth-order bi-Laplace equation. Since this discussion is mainly for theoretical purposes, we assume the simple harmonic extension followed by a decomposition to enforce $d\omega_{l+p}=0$ instead of a biharmonic extension. In Fig.~\ref{fgr:normalTangentialVF} (a), we illustrate the implementation of boundary conditions for the extension of normal harmonic forms to the interior cavity. In this evolving process, the outside surface is fixed and the inner cavity shrinks to null in order that the manifold with a cavity extends into a solid ball. Under the boundary condition Eq.~(\ref{eqn:bc_ext_n_1}) on the interior surface, the input normal harmonic forms (thin lines) are extended into the cavity, which also preserve curl-free properties shown as thick lines in Fig.~\ref{fgr:normalTangentialVF} (a).

Note that $d\mathfrak{E}(\omega)$ is a solution to the equation for solving  the extension of $d\omega$, by the uniqueness we impose, it must be $\mathfrak{E}(d\omega)$.  Thus, we can construct the following commutative diagram on the de Rham complexes for normal forms on the filtration of $M$:
\begin{center}
	\begin{tikzcd}
		\Omega^0_n(M_{0}) \arrow[d, "\mathfrak{E}_{0,1}"]
		\arrow[r, shift left, "d^0"] &
		\Omega^1_n(M_{0}) \arrow[d, "\mathfrak{E}_{0,1}"]
		\arrow[r, shift left, "d^1"]  &
		\Omega^2_n(M_{0}) \arrow[d, "\mathfrak{E}_{0,1}"]
		\arrow[r, shift left, "d^2"]  &
		\Omega^3_n(M_{0}) \arrow[d, "\mathfrak{E}_{0,1}"]  \\
		\Omega^0_n(M_{1}) \arrow[d, "\mathfrak{E}_{1,1}"] \arrow[r, shift left, "d^0"] &
		\Omega^1_n(M_{1}) \arrow[d, "\mathfrak{E}_{1,1}"]
		\arrow[r, shift left, "d^1"]  &
		\Omega^2_n(M_{1}) \arrow[d, "\mathfrak{E}_{1,1}"]
		\arrow[r, shift left, "d^2"]  &
		\Omega^3_n(M_{1}) \arrow[d, "\mathfrak{E}_{1,1}"]  \\
		\Omega^0_n(M_{2}) \arrow[d, "\mathfrak{E}_{2,1}"] \arrow[r, shift left, "d^0"] &
		\Omega^1_n(M_{2}) \arrow[d, "\mathfrak{E}_{2,1}"]
		\arrow[r, shift left, "d^1"]  &
		\Omega^2_n(M_{2}) \arrow[d, "\mathfrak{E}_{2,1}"]
		\arrow[r, shift left, "d^2"]  &
		\Omega^3_n(M_{2}) \arrow[d, "\mathfrak{E}_{2,1}"]  \\
		\cdots &
		\cdots &
		\cdots &
		\cdots
	\end{tikzcd}
\end{center}
which places the de Rham complex in the horizontal direction and the filtration-induced extensions in the vertical direction.

Now, we can discuss the direct relation of bases of normal harmonic forms induced by $\mathfrak{E}$. First, $\omega_n \in \ker d_l$ implies $\mathfrak{E}_{l,p}(\omega_n) \in \ker d_{l+p}$. Thus, there is an injective homomorphism from $\ker d_l$ to $\ker d_{l+p}$. This induces a homomorphism from the cohomology group $\ker d^k_l/ \im d^{k-1}_l$ to $\ker d^k_{l+p}/ \im d^{k-1}_{l+p}$, which, through de-Rham isomorphism between cohomology and harmonic spaces in $M_l$ and $M_{l_p}$, is equivalent to a homomorphism from the harmonic space $\mathcal{H}^k_{\Delta_n, l}$ to $\mathcal{H}^k_{\Delta_n, l+p}$. Instead of using the mapping between the equivalence classes, we can actually directly pick the unique harmonic representative $h_{n}\in \ker d^k \cup \ker \delta^{k+1}=\mathcal{H}^k_{\Delta_n}$ for each equivalence class in the cohomology, as we can pick the closed form that is orthogonal to $\im d^{k-1}$ which is $\ker \delta^{k}$ due to the adjointness between $d$ and $\delta$. However, for $h_n\in\mathcal{H}^k_{\Delta_n, l}$, its extension $\mathfrak{E}_{l,p}(h_n)$ is not necessarily an element of $\mathcal{H}^k_{\Delta_n, l+p}$. Nevertheless, composed with the simple $L_2$ projection onto the finite dimensional normal harmonic space $P_{\mathcal{H}^k_{\Delta_n, l+p}}$, we have the linear map (also a homomorphism) $\Psi_{n,l,p}=P_{\mathcal{H}^k_{\Delta_n, l+p}}\circ\mathfrak{E}_{l,p}:\mathcal{H}^k_{\Delta_n, l}\to \mathcal{H}^k_{\Delta_n, l+p}$.

The map between these two normal harmonic spaces is neither necessarily injective nor necessarily surjective. In fact, if  $h_n \in \mathcal{H}^k_{\Delta_n,l}$ is not in $\im \Psi_{n,l-1,1}$, it is said to be born at index $l$; if $p$ is the smallest integer such that $\Psi_{n,l,p}(h_n)=0$, it is said to die at index $l+p$, with a persistence of $p$. This is consistent with the persistence of the relative cohomology $H^k(M,\partial M)$ and the (absolute) homology $H_{3-k}(M)$.

\subsubsection{Tangential harmonic forms}
As there is a one-to-one correspondence between tangential $k$-forms and normal $(3\!-\!k)$-forms, it is indeed sufficient to study the tangential forms only. For completeness and flexibility in numerical implementation, we provide a brief discussion on this dual case.


We first note that there is a homomorphism from coclosed forms on $M_l$ to coclosed forms on $M_{l+p}$, i.e., from $\ker \delta_{l}$ to $\ker \delta_{l+p}$ when restricted to tangential forms $\Omega_t(M_l)$. The same harmonic extension $\mathfrak{E}_{l,p}$ can be obtained through the minimization of the Dirichlet energy $\langle d\omega,d\omega\rangle + \langle \delta \omega,\delta \omega\rangle$ in $M^c_{l,p}$. For tangential forms, $\star\omega_l|_{\partial M_l}=0$, we also impose the condition that $d \omega$ is tangential to the boundary ($\star d \omega_l|_{\partial M_l}=0$). We keep this condition on $\partial M_{l+p}$, on $\partial M_l$ we impose continuity $\omega_{l+p}|_{\partial M_l} = \omega_l|_{\partial M_l}$ and $ d\omega_{l+p}|_{\partial M_l} = d \omega_l|_{\partial M_l}$. A unique solution is again found by forcing it to have 0 projection to the kernel of a mixed-type boundary condition Laplace equation~\cite{zhao20193d}. 

To illustrate it with a tangential 1-form $\omega_l$, we can impose the homogeneous boundary condition for the proxy vector field $\bv$ on $\partial M_{l+p}$ as in Eq.~(\ref{eqn:bdyT1form}),
\begin{equation}
\label{eqn:bc_ext_t_1}
\bv_{l+p}\cdot\bn= 0 , \quad \nabla_{\bn}(\bv_{l+p}\cdot\bt_1)=0,\quad \nabla_{\bn}(\bv_{l+p}\cdot\bt_2)=0;
\end{equation}
whereas on $\partial M_l$, the Dirichlet boundary condition $\bv_{l+p}=\bv_l$ is equivalent to
\begin{equation}
\label{eqn:bc_ext_t}
\bv_{l+p}\cdot\bt_1=\bv_{l}\cdot\bt_1,\quad \bv_{l+p}\cdot\bt_2=\bv_{l}\cdot\bt_2,\quad \bv_{l+p}\cdot\bn=0.
\end{equation}
In this case, we can enforce $\mathfrak{E}_{l,p}(\ker \delta_l)\subset \ker \delta_{l+p}$. For example, Fig.~\ref{fgr:normalTangentialVF} (b) shows the extension of tangential harmonic forms from a torus to a solid sphere where both boundary conditions Eqs.~(\ref{eqn:bc_ext_t_1}) and (\ref{eqn:bc_ext_t}) are applied. The inputs (thick lines) are only circulations shown in the right chart of Fig.~\ref{fgr:normalTangentialVF} (b), while the extended outputs (thin lines) are tangential harmonic forms as well. Therefore, we can construct the following commutative diagram on the de Rham complexes for tangential forms on the filtration of $M$:
\begin{center}
	\begin{tikzcd}
		\Omega^0_t(M_{0}) \arrow[d, "\mathfrak{E}_{0,1}"]
		&
		\Omega^1_t(M_{0}) \arrow[d, "\mathfrak{E}_{0,1}"]
		\arrow[l, shift left, "\delta^1"]  &
		\Omega^2_t(M_{0}) \arrow[d, "\mathfrak{E}_{0,1}"]
		\arrow[l, shift left, "\delta^2"]  &
		\Omega^3_t(M_{0}) \arrow[d, "\mathfrak{E}_{0,1}"]
		\arrow[l, shift left, "\delta^3"]  \\
		\Omega^0_t(M_{1}) \arrow[d, "\mathfrak{E}_{1,1}"] &
		\Omega^1_t(M_{1}) \arrow[d, "\mathfrak{E}_{1,1}"]
		\arrow[l, shift left, "\delta^1"]  &
		\Omega^2_t(M_{1}) \arrow[d, "\mathfrak{E}_{1,1}"]
		\arrow[l, shift left, "\delta^2"]  &
		\Omega^3_t(M_{1}) \arrow[d, "\mathfrak{E}_{1,1}"] \arrow[l, shift left, "\delta^3"]  \\
		\Omega^0_t(M_{2}) \arrow[d, "\mathfrak{E}_{2,1}"]  &
		\Omega^1_t(M_{2}) \arrow[d, "\mathfrak{E}_{2,1}"]
		\arrow[l, shift left, "\delta^1"]  &
		\Omega^2_t(M_{2}) \arrow[d, "\mathfrak{E}_{2,1}"]
		\arrow[l, shift left, "\delta^2"]  &
		\Omega^3_t(M_{2}) \arrow[d, "\mathfrak{E}_{2,1}"] \arrow[l, shift left, "\delta^3"] \\
		\cdots &
		\cdots &
		\cdots &
		\cdots
	\end{tikzcd}
\end{center}

Similar to the normal form case, through the composition with the simple $L_2$ projection onto the finite dimensional tangential harmonic space $P_{\mathcal{H}^k_{\Delta_t, l+p}}$, we have a linear map (also a homomorphism) between the tangential harmonic spaces of different manifolds in the filtration, $\Psi_{t,l,p}=P_{\mathcal{H}^k_{\Delta_t, l+p}}\circ\mathfrak{E}_{l,p}:\mathcal{H}^k_{\Delta_t, l}\to \mathcal{H}^k_{\Delta_t, l+p}$. If  $h_t \in \mathcal{H}^k_{\Delta_t,l}$ is not in $\im \Psi_{t,l-1,1}$, it is said to be born at index $l$. If $p$ is the smallest integer such that $\Psi_{t,l,p}(h_t)=0$, it is said to die at index $l+p$, with a persistence of $p$. This is consistent with the persistence of the (absolute) cohomology $H^k(M)$ and the relative homology $H_{3-k}(M,\partial M)$.

\subsubsection{Relation among persistent cohomologies under different boundary conditions}

As discussed in section~\ref{Discreteformandspectrumanalysis}, with the duality through Hodge star,
there are only three independent singular spectra $T$, $N$ and $C$ for the three differential/codifferential operators (two for gradient operators under tangential or normal conditions, and one curl operator with either tangential or normal boundary condition). The unions of these spectra produce all the eigenvalues of the eight possible Hodge Laplacians on an arbitrary compact manifold $M$ embedded in a flat 3D space. Moreover, the intersections of spaces spanned by left or right singular vectors of singular value 0 for these operators form the tangential and normal harmonic spaces. Thus, we can restrict our discussion to either normal or tangential fields without loss of generality.

We now discuss the persistence from the perspective of evolving Hodge Laplacian operators. Note that the following discussion is to provide theoretical backgrounds for our proposed use of the evolution of eigenvalues, but not for implementations, since some of the operators discussed may not be sparse matrices when discretized. Recall that for any two manifolds $M_l$ and $M_{l+p}$ in any type of filtration, there is an inclusion map $\mathfrak{I}_{l,p}:M_l \hookrightarrow M_{l+p}$. We call $M_{l+p}$ the \emph{$p$-evolution manifold} of $M_{l}$. We can directly investigate whether a harmonic form in $M_l$ survived in its $p$-evolution manifold, by defining a restricted subset $\tilde{\Omega}^{k}_{p}(M_l)$ of $\Omega^k(M_{l+p})$ and using it to define modified differential and codifferential operators on $M_l$. This restricted subset is given by
\begin{equation}
\label{eqn:strictionSet}
\tilde{\Omega}^{k}_{p}(M_l)=\{\omega\in\Omega^{k}(M_{l+p}) | d^{k}_{l+p}\omega \in \mathfrak{E}_{l,p}(\ker d^{k+1}_{l}) \}.
\end{equation}
This space can be equipped with a modified operator $\tilde{d}^{k}_{l+p}$ that maps it to $\Omega^{k+1}(M_l)$, which is defined as the compound of $d^{k}_{l+p}$ followed by the pullback through the inclusion, i.e., $\tilde{d}^k_{l+p}=\mathfrak{I}_{l,p}^*\circ d^k_{l+p}$. Assuming that we use normal differential forms, we have $d^{k+1}_{l}\tilde{d}^{k}_{l+p}=0$ on $\tilde{\Omega}^k_p(M_l)$ as a result of the definition of the restricted space.
For $\omega\in \Omega^{k-1}(M_l)$, we have $d^{k-1}_{l+p} \mathfrak{E}_{l,p}(\omega)=\mathfrak{E}_{l,p}(d^{k-1}_{l+p}\omega)\in \mathfrak{E}_{l,p}(\ker d^k_l)$, thus $\mathfrak{E}_{l,p}(\Omega^{k-1}(M_l)) \subseteq \tilde{\Omega}^{k-1}_{p}(M_l)$ for $p\geq 0$. Therefore, we can construct the following the $p$-evolution differential form diagram
\begin{center}
	\begin{tikzcd}
		\Omega^0(M_{l  }) \arrow[r, shift left, "d^0_l"] \arrow[dd,"\mathfrak{E}_{l,p}"] &
		\Omega^1(M_{l  }) \arrow[ldd, shift left, "\tilde{\delta}^1_{l+p}"]
		\arrow[r, shift left, "d^1"] \arrow[l, shift left, "\delta^1_l"] \arrow[dd,"\mathfrak{E}_{l,p}"] &
		\Omega^2(M_{l  }) \arrow[ldd, shift left, "\tilde{\delta}^2_{l+p}"]
		\arrow[r, shift left, "d^2_l"] \arrow[l, shift left, "\delta^2_l"] \arrow[dd,"\mathfrak{E}_{l,p}"] &
		\Omega^3(M_{l  }) \arrow[ldd, shift left, "\tilde{\delta}^3_{l+p}"]
		\arrow[l, shift left, "\delta^3_l"] \\
		& & & \\
		\tilde{\Omega}^{0}_{p}(M_l) \arrow[ruu, shift left, "\tilde{d}^0_{l+p}"] &
		\tilde{\Omega}^{1}_{p}(M_l) \arrow[ruu, shift left, "\tilde{d}^1_{l+p}"] &
		\tilde{\Omega}^{2}_{p}(M_l) \arrow[ruu, shift left, "\tilde{d}^2_{l+p}"] &
	\end{tikzcd}
\end{center}
where $\tilde{\delta}^k_{l+p}$ denotes the adjoint operator of $\tilde{d}^k_{l+p}$. Based on this diagram,
the $p$-evolution Hodge Laplacian $\Delta^{k}_{l,p}$: $\Omega^k(M_l)\rightarrow\Omega^k(M_l)$ can be defined on $M_l$ as
\begin{equation}
\label{eqn:restrictLaplacian}
\Delta^{k}_{l,p} = \delta^{k+1}_{l} d^k_l+\tilde{d}^{k-1}_{l+p}\tilde{\delta}^{k}_{l+p},
\end{equation}
which leads to the definition of the $p$-evolution harmonic space as $\mathcal{H}^k_{l,p}=\ker{\Delta^{k}_{l,p}}=\ker{d^k_l}\cap\ker{\tilde{\delta}^{k}_{l+p}}$. The $p$-evolution (tangential) $k$-form spectra are the sets of $\Delta^{k}_{l,p}$'s eigenvalues for $k=0,1,2,3$. By comparing the $p$-evolution Laplace operator $\Delta^{k}_{l,p}$
and the Laplace operator $\Delta^{k}_{l,0}$, the eigenvalues of the unmodified part, $\delta^{k+1}_{l} d^k_l$, are preserved, and the eigenvalues involving the pullback of the restricted operators are varying with $p$. Next, we examine the part involving $\tilde{d}^{k-1}_{l+p}\tilde{\delta}^{k}_{l+p}$. For any $\alpha\in \ker \tilde{\delta}^k_{l+p},$ and any $\tilde\beta\in \tilde{\Omega }_p^{k-1}(M_l)$, we have $0=\langle\tilde{\delta}^k_{l+p}\alpha,\tilde{\beta}\rangle = \langle\alpha,\tilde{d}^{k-1}_{l+p}\tilde\beta\rangle$. For any $\beta \in\Omega^{k-1}(M_l)$, we have $\langle\delta^k_l \alpha,\beta\rangle = \langle\alpha,d^{k-1}_{l}\beta\rangle = \langle\alpha,\tilde{d}^{k-1}_{l+p}\mathfrak{E}_{l,p}(\beta)\rangle=0$. Therefore, $\ker\tilde{\delta}^k_{l+p}\subset\ker\delta^k_l\subset\Omega^k(M_l)$.

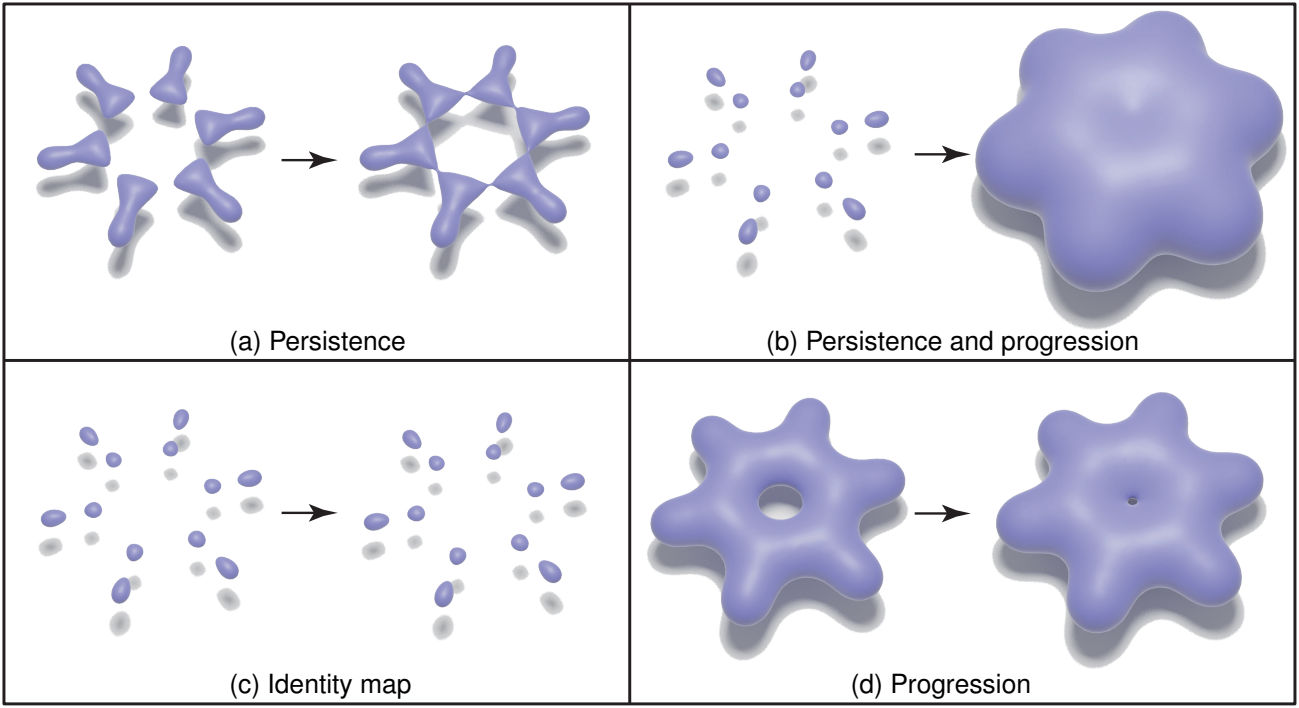
\begin{figure*}[t!]
	\centering
	\includegraphics[height=3.7in]{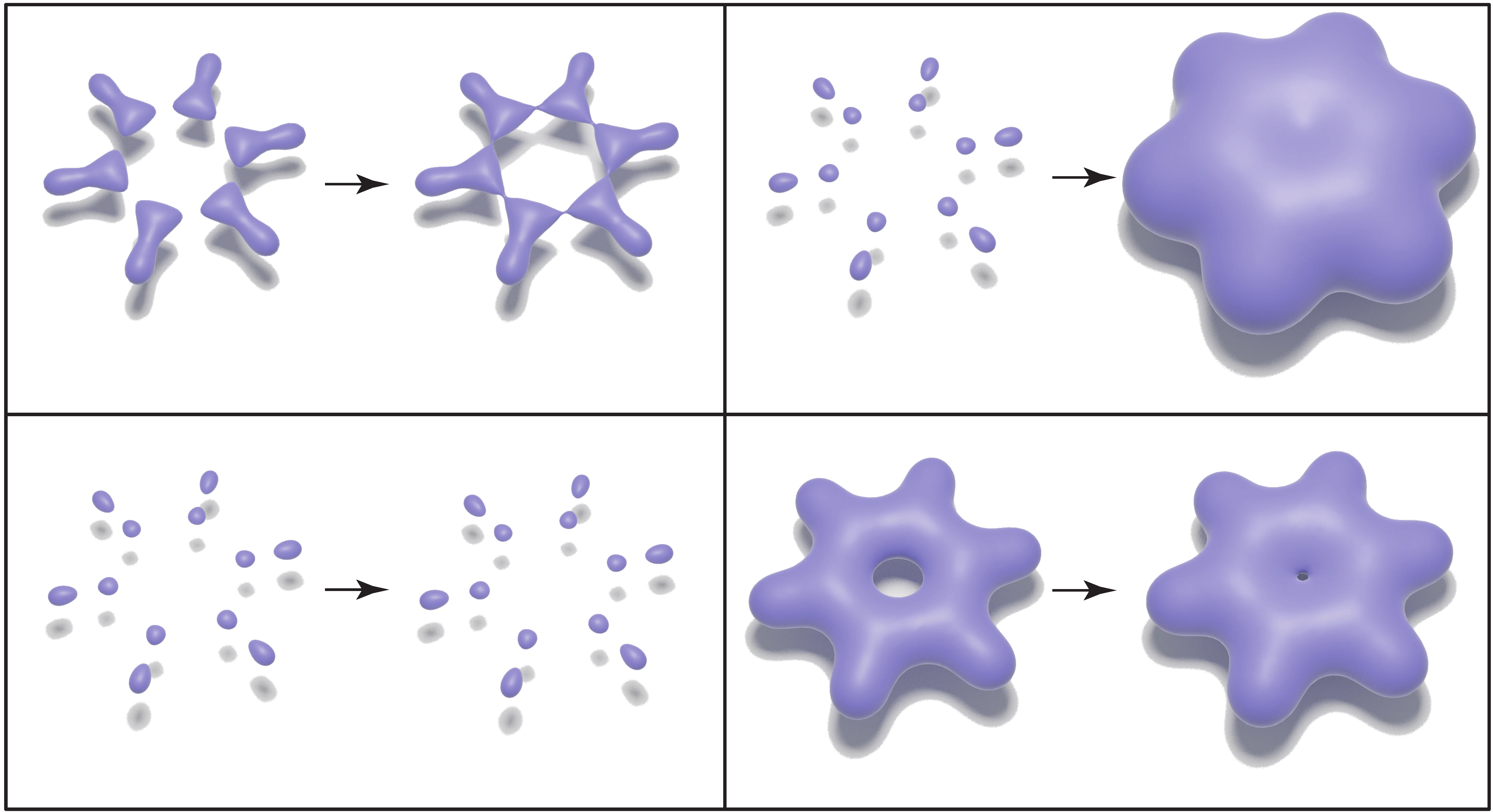}
	\begin{picture}(17,0)
	\put( -150, 147){\text{(a) Persistence}}
	\put(   53, 147){\text{(b) Persistence and progression}}
	\put( -150, 17){\text{(c) Identity map}}
	\put(   85, 17){\text{(d) Progression}}
	\end{picture}
	\caption{Persistence and progression on benzene.}
	\label{fgr:persistenceProgression}
\end{figure*}

Thus, in terms of persistent cohomology, we may examine the kernel of $p$-evolution Laplace operator for the persistence of topological features of $M_l$ in $M_{l+p}$. In the perspective of spectral analysis, this change is reflected in the multiplicity of the eigenvalue $0$, which changes if $\dim{(\ker\tilde{\delta}^k_l)}<\dim{(\ker\delta^k_l)}$, or remains unchanged when $\dim{(\ker\tilde{\delta}^k_l)}=\dim{(\ker\delta^k_l)}$. In the former case, as shown in Fig.~\ref{fgr:persistenceProgression} (a),  multiplicity of 0 (the number of connected components) is reduced for $\Delta^0_{l,p}$, whereas $\Delta^1_{l,p}$ has a new 0 (a tunnel) that is not present in $\Delta^1_{l,p}$. For the latter case, the inclusion map is homotopic to a geometrical deformation of the manifold, which implies the same topology. Fig.~\ref{fgr:persistenceProgression} (d) illustrate an example where the size of tunnel shrinks, and the cohomology groups are isomorphic.

The spectra are continuous when corresponding manifolds are continuously deforming, since, as discussed above, when the level set values are close, the deformation is close to an isometric, and the eigenvalues of Hodge Laplacian is determined by the metric tensor. In particular, the smallest non-zero eigenvalues are continuous if the dimension of null space is stable, but are typically non-differentiable when the multiplicity of eigenvalue 0 is changed. The birth of non-zero eigenvalues is the death of topological features, which signals the death of harmonic basis fields; whereas the birth of zero eigenvalues indicates the birth of topological features. Moreover, the changes in leading smallest non-zero eigenvalues can thus indicate possible pending topological changes as well as the geometric properties when the manifold evolves without topological changes.

For instance, for the $l$-th manifold of the filtration of $M$, $\{\lambda^{T}_{l,i}\}$, $\{\lambda^{C}_{l,i}\}$ and $\{\lambda^{N}_{l,i}\}$ give the eigenvalues of the $T$, $C$ and $N$ sets respectively. In particular, the multiplicities of the zero eigenvalues in   $\lambda^{T}_{l,0}$, $\lambda^{C}_{l,0}$, and $\lambda^{N}_{l,0}$ are associated with Betti numbers $\beta_0, \beta_1$ and $\beta_2$, respectively. Additionally,  $\lambda^{T}_{l,1}$, $\lambda^{C}_{l,1}$, and $\lambda^{N}_{l,1}$ are the first non-zero eigenvalues, which are known as the Fiedler values in graph theory, an indicator of how well the graph is connected. 

In summary, the correspondence established by the spectral analysis provides us with tools to investigate both types of manifold evolution, with persistence for topological features and spectral progression for the geometric properties.

\section{Evolutionary de Rham-Hodge analysis of geometric shapes}\label{resDHT}
 
 In this section, we present the application of the proposed evolutionary de Rham-Hodge method.  We demonstrate the spectral analysis with evolutionary de Rham Laplace operators and illustrate their topological persistence and geometric progression associated with submanifolds in $\mathbb{R}^3$. The evolving manifolds in our studies are generated by applying Eq.~(\ref{eqn:manifolds}) to point cloud datasets with a varying level set $c$, with a fixed scaling parameter $\eta$. 

For clarity, the first three examples are simple point sets  consisting of few points. The two-body set has the location coordinates in $\{(-1.5, 0, 0), (1.5, 0, 0)\}$, and for the four-body and eight-body sets. We duplicate the two-body set by translating $\pm 1.5$ along the $y$-axis, and duplicate the four-body set by translating $\pm 1.5$ along the $z$-axis respectively. Next, we present two concrete molecular examples with interesting topological and geometric features, benzene ($\text{C}_6\text{H}_6$) and fullerene ($\text{C}_{60}$). We show in these proof-of-concept examples that the evolution of leading smallest eigenvalues provides additional information to that of the persistent Betti numbers, which are the same as those of persistent homology analysis. That is, we propose to extend the evaluation of the manifold evolution from persistent Betti numbers (i.e., the multiplicity of the zero eigenvalues of evolutionary  de Rham Laplace operators) to a larger subset of the spectra.

\subsection{Two-body system}

\begin{figure}
	\centering
	\includegraphics[height=1.7in]{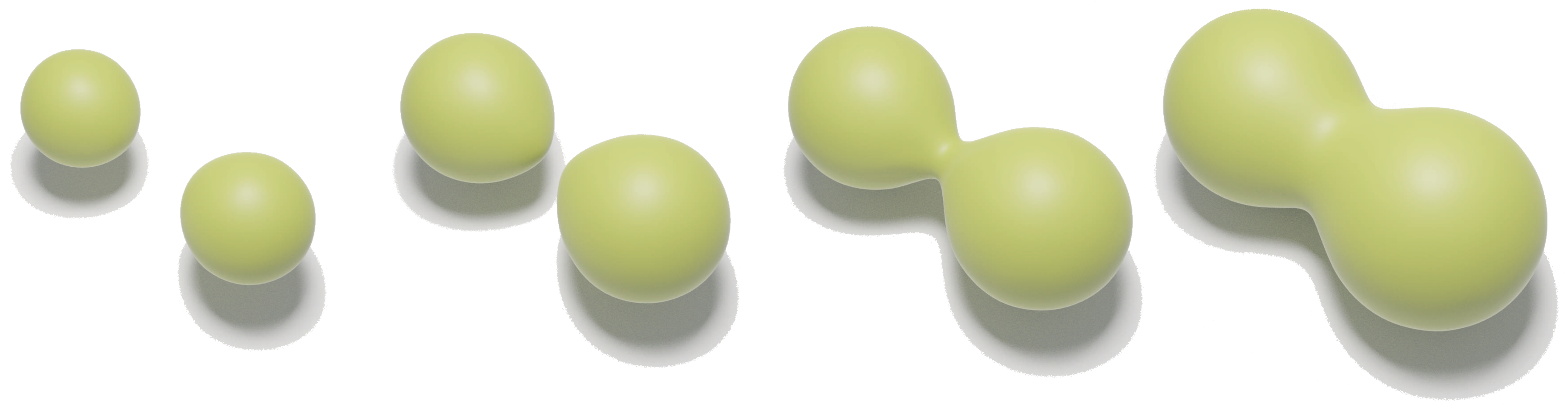}
	\begin{picture}(0,0)
	\put( -485, 110){\textbf{a}}
	\put( -370, 110){\textbf{b}}
	\put( -250, 110){\textbf{c}}
	\put( -135, 110){\textbf{d}}
	\end{picture}
	\caption{Snapshots of evolving manifold with the two-body system. {\bf a}, {\bf b}, {\bf c} and {\bf d} are snapshots from the beginning to the end. {\bf b} and {\bf c} show the transition of the Betti-0 number from 2 to 1.
	}
	\label{fgr:twoAtoms}
\end{figure}
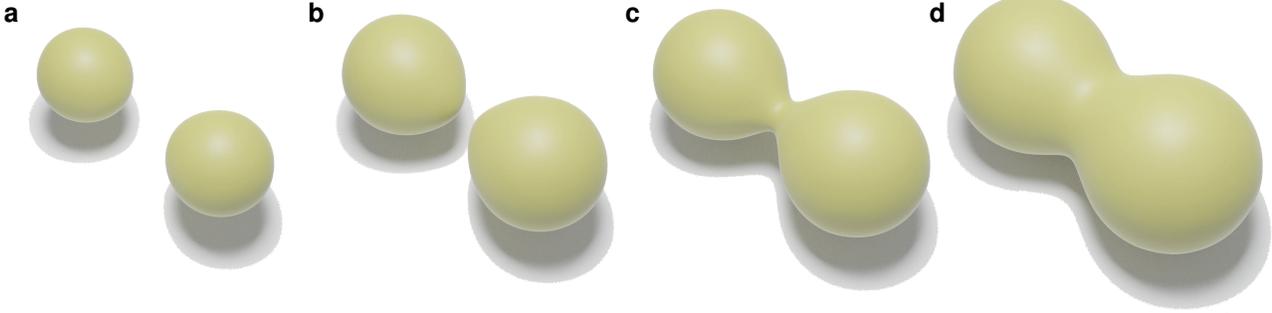

\begin{figure}
	\centering
	\includegraphics[height=1.82in]{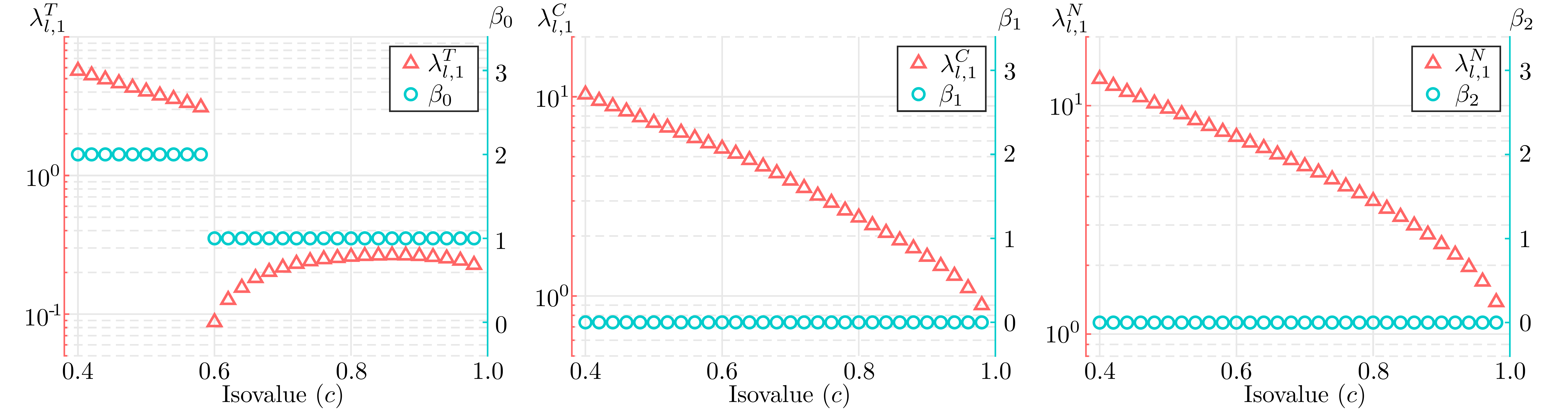}
	\begin{picture}(0,0)
	\put(-160, 135){\textbf{i}}
	\put(   5, 135){\textbf{ii}}
	\put( 165, 135){\textbf{iii}}
	\end{picture}
	\caption{ Eigenvalues and Betti numbers  vs isovalue ($c$) of the two-body system with $\eta = 1.19$ and $\max(\rho)\approx 1.0$. 
	{\bf i} shows the smallest eigenvalues of the $T$ set. 
	The drops at $c=0.6$ correspond to snapshots in Figs.~\ref{fgr:twoAtoms} {\bf b} and {\bf c}. 
	{\bf ii} and {\bf iii} show the smallest eigenvalues of the $C$ and $N$ sets respectively.
	}
	\label{fgr:twoAtomsEigBetti}
\end{figure}
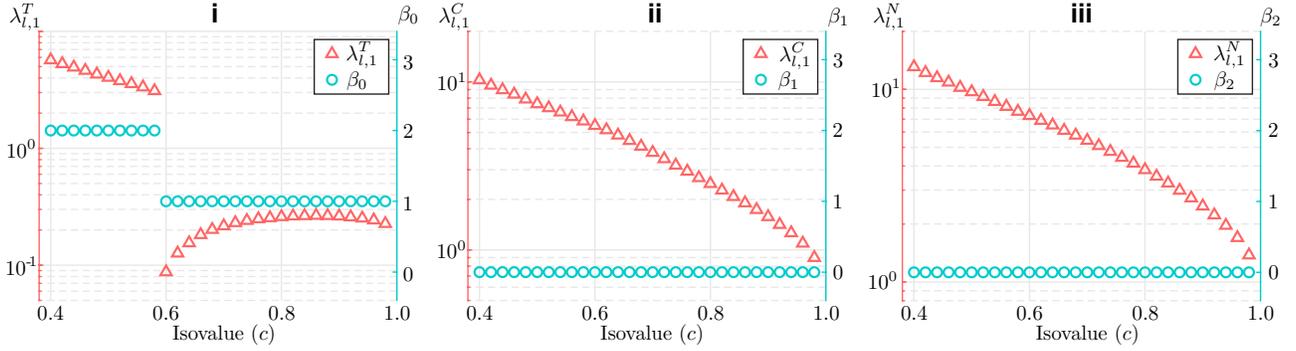

Our first example illustrates the evolving manifold with a two-body system, in which the initial two connected components merge into one. In this evolution, only the number of components persistent $\beta_0$ changes from 2 to 1, with the other Bettie numbers remain at $0$ throughout. As shown in Fig.~\ref{fgr:twoAtoms}, the two connected components gradually approach each other as the isovalue grows and eventually touch each other as more volume is enclosed.

The change in topology can be observed directly from the blue circle plots in Fig.~\ref{fgr:twoAtoms}, where persistent $\beta_0$ is dropped from 2 to 1 when $c$ increased to around 0.6, and the curves for persistent $\beta_1$ and $\beta_2$ remained flat due to the lack of tunnels or cavities in the system. However, the persistent Betti numbers do not provide any information about the volume increase of the manifold during the evolution, or the increase in the size of the tube-like structure between the two blobs around the body centers after they touch. In contrast, the orange triangles in Fig.~\ref{fgr:twoAtomsEigBetti} show how the first nonzero eigenvalues (Fiedler values) in the three singular spectra ($T$, $C$ and $N$) demonstrated both the topological transition and geometric progression in the evolving manifold.

First, one may observe that the discontinuity for the Fiedler values of the tangential gradient fields $T$ coincides with the jump of persistent $\beta_0$ in Fig.~\ref{fgr:twoAtomsEigBetti} {\bf i}, whereas the Fiedler values of the tangential/normal curl fields $C$ and that of the normal gradient fields $N$ are both smooth as shown in Figs.~\ref{fgr:twoAtomsEigBetti} {\bf ii} and {\bf iii}. These behaviors are consistent with the evolution process only having changes in the number of connected components. More precisely, the multiplicity of the eigenvalue zero in $T$ is $\beta_0=2$ at the beginning, so the Fiedler values can be seen as the third eigenvalue, whereas after the merging, it is switched to be the second eigenvalue, which contributes to the discontinuity in its value. As we will see in later examples, this behavior for the persistence to be directly observable in the discontinuity of Fiedler values happening at the same isovalue when the Betti numbers jump to different integers is generic, which indicates that the birth of non-zero eigenvalue and the death of the harmonic basis are both linked to the death of topological features (homology generators). Moreover, as the tube between the two blobs is created, the extreme values of the first oscillation mode can be placed further apart along the line connecting the two atoms. Thus, $\lambda^T_{l,1}$ jumps to a small value. It grows as the structure becomes stiffer when the narrow tube turns thicker before it eventually decays again as the entire shape turns softer as a ball with a growing radius. Figs.~\ref{fgr:twoAtomsEigBetti} {\bf ii} and {\bf iii} show the smoothness of $\lambda^C_{l,1}$ and $\lambda^N_{l,1}$ which is consistent with the invariant 1st and 2nd Betti numbers.

\subsection{Four-body system}

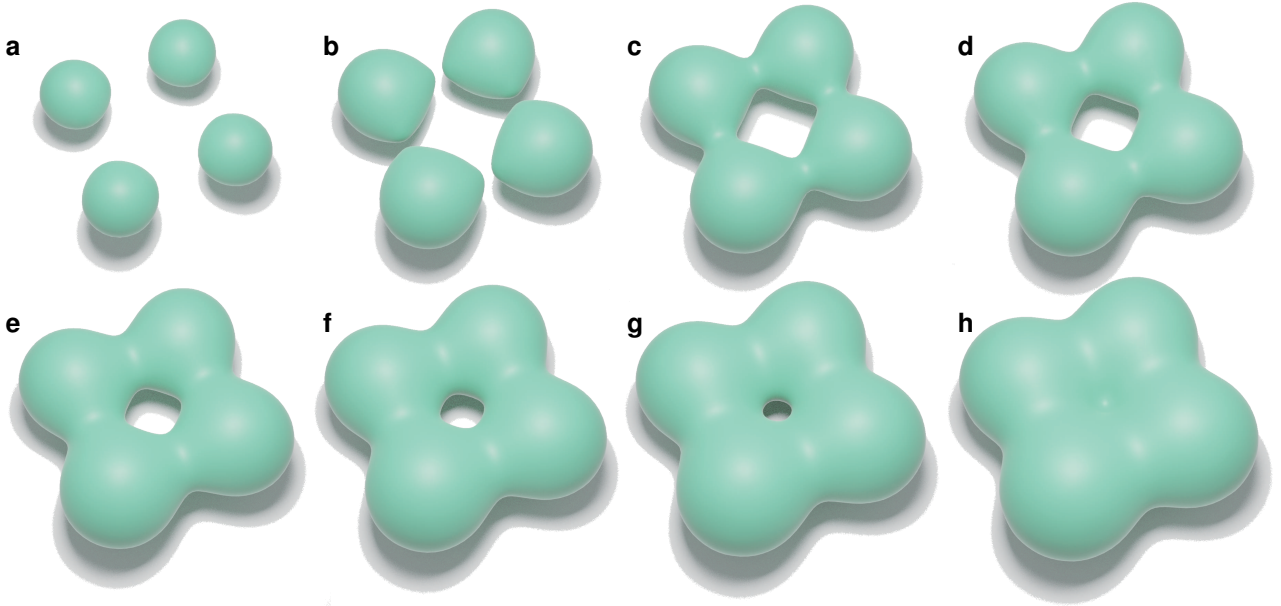
\begin{figure}
	\centering
	\includegraphics[height=3.2in]{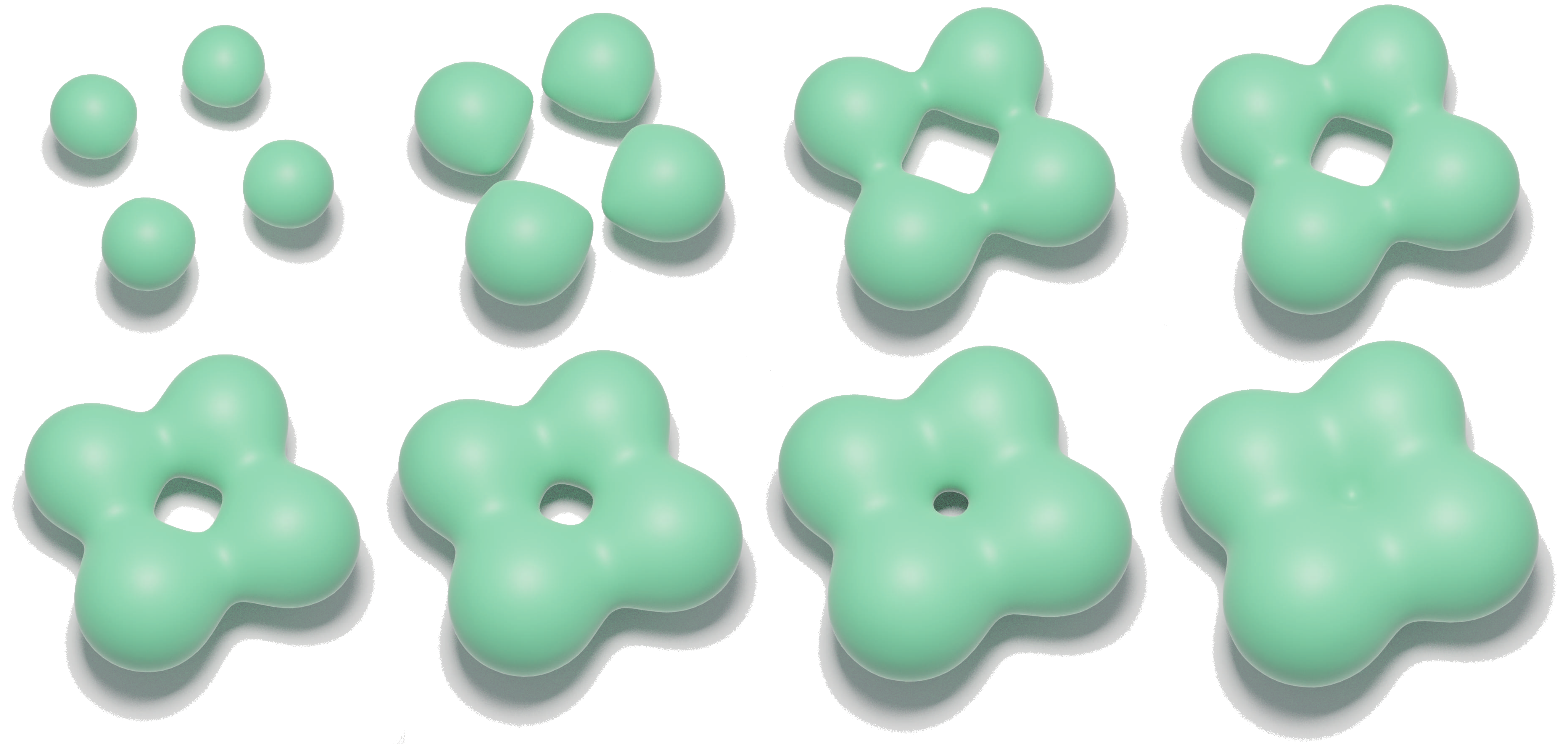}
	\caption{ Snapshots of evolving manifolds with the four-body system. {\bf a} is the initial point of four components; {\bf b} and {\bf c} show the transition of a ring formed and the persistent Betti-0 number changes from 4 to 1. {\bf g} and {\bf h} show the vanishing of the ring and the persistent Betti-1 number changes from 1 to 0.
	}
	\label{fgr:fourAtoms}
	\begin{picture}(0,0)
	\put( -240, 255){\textbf{a}}
	\put( -120, 255){\textbf{b}}
	\put(   -5, 255){\textbf{c}}
	\put(  120, 255){\textbf{d}}
	\put( -240, 150){\textbf{e}}
	\put( -120, 150){\textbf{f}}
	\put(   -5, 150){\textbf{g}}
	\put(  120, 150){\textbf{h}}
	\end{picture}
\end{figure}

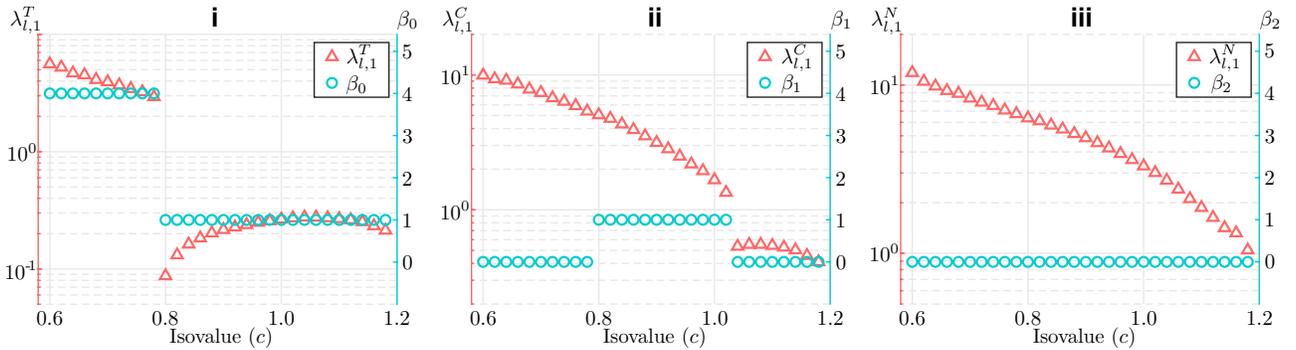
\begin{figure}
	\centering
	\includegraphics[height=1.82in]{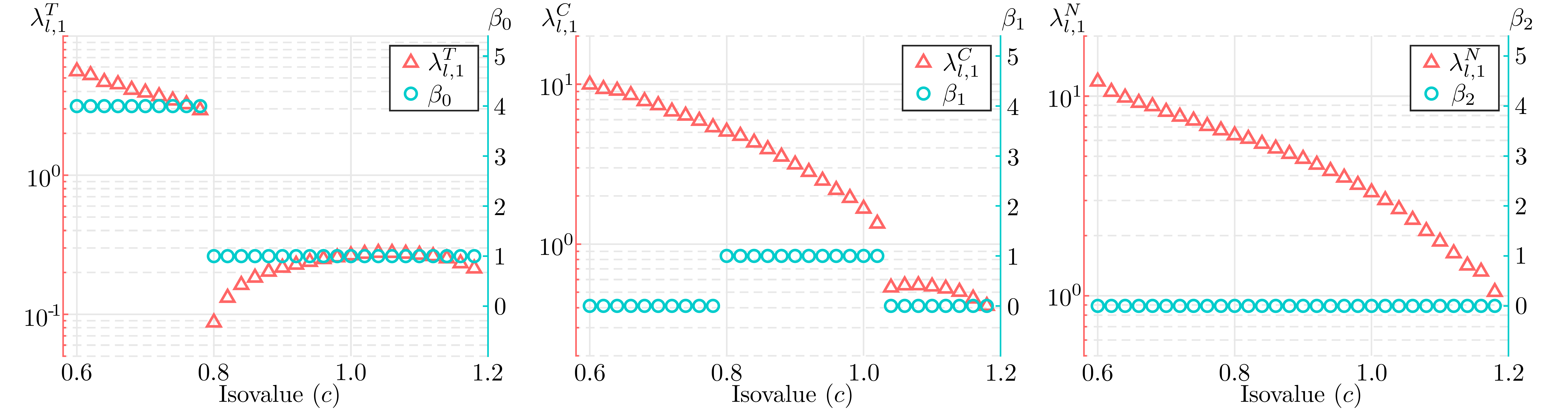}
	\begin{picture}(0,0)
	\put(-160, 135){\textbf{i}}
	\put(   5, 135){\textbf{ii}}
	\put( 165, 135){\textbf{iii}}
	\end{picture}
	\caption{Eigenvalues and Betti numbers  vs isovalue ($c$) of the four-body system with $\eta = 1.19$ and $\max(\rho)\approx 1.2$. 
	{\bf i} shows the smallest eigenvalues of the $T$ set. At near $c=0.80$, the persistent Betti-0 number changes from 4 to 1.
	{\bf ii} shows the smallest eigenvalues of the $C$ set. At around $c=1.02$, the persistent Betti-1 number changes from 1 to 0. 
	{\bf iii} shows the smallest eigenvalues of the $N$ set.
	}
	\label{fgr:fourAtomsEigBetti}
\end{figure}

As another example, we explore an evolution that involves changes in both the number of components persistent $\beta_0$ and the number of tunnels $\beta_1$. With two points added to the two-body set to form a planar square, the evolving manifold can contain a tunnel for a range of isovalues, when each of the four components touches two neighbors to form a ring, which will eventually disappear as the level set value increases to the point that the tunnel in the middle is filled. During the same process, persistent $\beta_0$ drops from four to one when persistent $\beta_1$ increased to one with the formation of the tunnel, but persistent $\beta_0$ stays at $1$ when persistent $\beta_1$ changes back to zero with the disappearance of the tunnel. The  persistent Betti number $\beta_2$ remains unchanged as there is no cavity in the system.

In terms of the geometric measurements, the total volume continuously increases, and once the tunnel appears, the size of the handle dual to the tunnel also increases. Finally, at the time of disappearing of the tunnel, two concave surfaces are formed on each side of the blocked tunnel with the concavity decreases with an increasing level set parameter.

 Fig.~\ref{fgr:fourAtomsEigBetti} shows all the Fiedler values varying over time, along with the relevant Betti numbers.
As both $\beta_0$ and $\beta_1$ change during the evolution, $\lambda^T_{l,1}$ and $\lambda^C_{l,1}$ are non-differentiable for this example. On the other hand, $\beta_2$ is invariant and thus $\lambda^N_{l,1}$ is smooth. Fig.~\ref{fgr:fourAtomsEigBetti} {\bf i} exhibits a similar pattern as the two-body case of  $\lambda^T_{l,1}$. As the volume of the manifold increases, $\lambda^T_{l,1}$ decays until the four components are connected, at which point $\lambda^T_{l,1}$ drops to a much smaller value. After the discontinuity, the increasing handle size leads to an initial growth of $\lambda^T_{l,1}$ due to the increased stiffness of the system, before returning to the decreasing trend as the system becomes more flexible with the increase in the overall volume. In Fig.~\ref{fgr:fourAtomsEigBetti} {\bf ii}, one may observe the difference compared with the first case as we introduce the changes in persistent $\beta_1$. When $\beta_1$ changes from zero to one through the connection of the four components, $\lambda^C_{l,1}$ does not actually change much, because the tangential/normal curl field is not largely influenced when the handle size is nearly zero. In stark contrast, $\lambda^C_{l,1}$ is discontinuous when $\beta_1$ changes back down to zero as the hole disappears. The behavior of $\lambda^C_{l,1}$ after the discontinuity is similar to that of $\lambda^T_{l,1}$, an initial increase in stiffness and then a decrease again. Moreover, by comparing Figs.~\ref{fgr:fourAtomsEigBetti} {\bf i} and {\bf ii}, we observe that the value of $\lambda^T_{l,1}$ starts to decrease just when $\lambda^C_{l,1}$ is discontinuous, as the structural change in the tunnel also contributed to the ``stiffness'' of the tangential gradients. Finally, Fig.~\ref{fgr:fourAtomsEigBetti} {\bf iii} shows the smooth Fiedler values $\lambda^N_{l,1}$ with an unchanged persistent $\beta_2$.

In summary, from the second example, one can notice that $\lambda^C_{l,1}$ can reveal the information of persistent $\beta_1$
and some geometric properties after the disappearance of the hole. In addition, the coincidental topological changes, the birth of hole that coincides with the death of a few connected components, can be distinguished by the spectral functions $\lambda^T_{l,1}$ and $\lambda^C_{l,1}$.

\subsection{Eight-body system}
\begin{figure}
	\centering
	\includegraphics[height=2.4in]{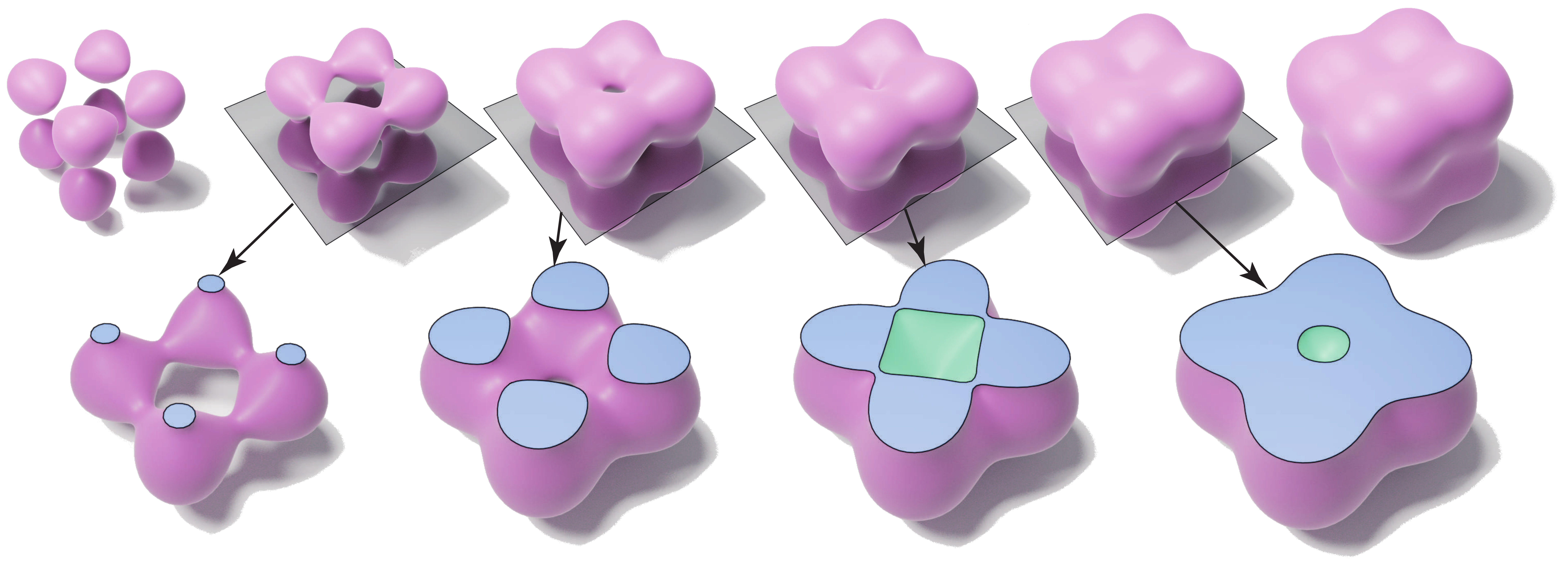}
	\begin{picture}(0,0)
	\put( -487, 157){\textbf{a}}
	\put( -410, 157){\textbf{b}}
	\put( -465, 83){\textbf{b}'}
	\put( -333, 157){\textbf{c}}
	\put( -357, 83){\textbf{c}'}
	\put( -255, 157){\textbf{d}}
	\put( -244, 83){\textbf{d}'}
	\put( -176, 157){\textbf{e}}
	\put( -125, 83){\textbf{e}'}
	\put( -98, 157){\textbf{f}}
	\end{picture}
	\caption{Snapshots of evolving manifold with the eight-body system. 
	{\bf a} presents the initial state with eight components.
    {\bf b} and {\bf c} show the formation of 6 tunnels when the persistent Betti-0 number changes from 8 to 1, and the persistent Betti-1 number changes from 0 to 5. 
     {\bf d} and {\bf e} illustrate that a cavity appears, so the persistent Betti-1 number drops to 0 and the persistent Betti-2 number increases to 1. 
     {\bf f} shows a solid volume without cavity. The gray planes cut manifolds to create  cross-section views to illustrate the process of the formation of cavity as shown in {\bf b'}, {\bf c'}, {\bf d'} and {\bf e'}.
	}
	\label{fgr:eightAtoms}
\end{figure}
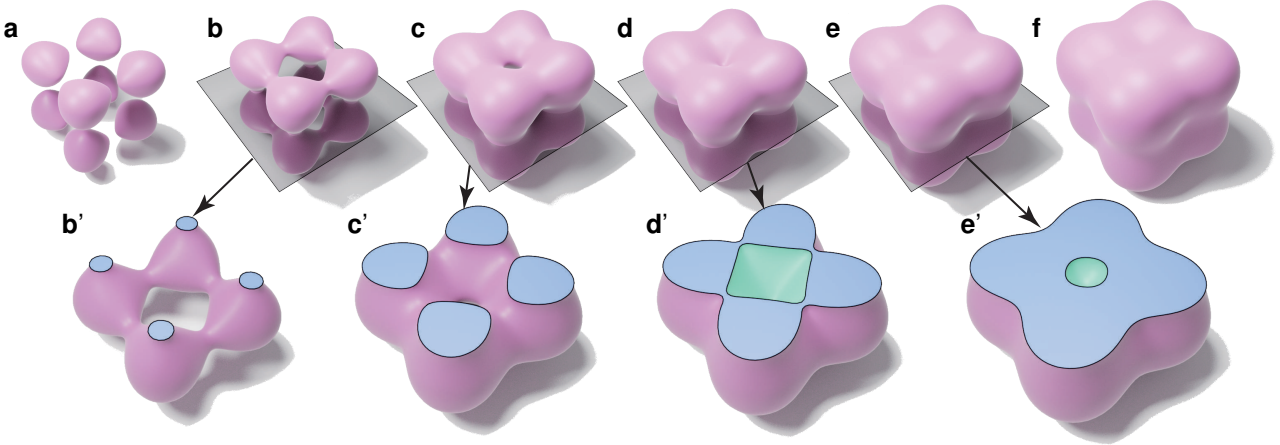
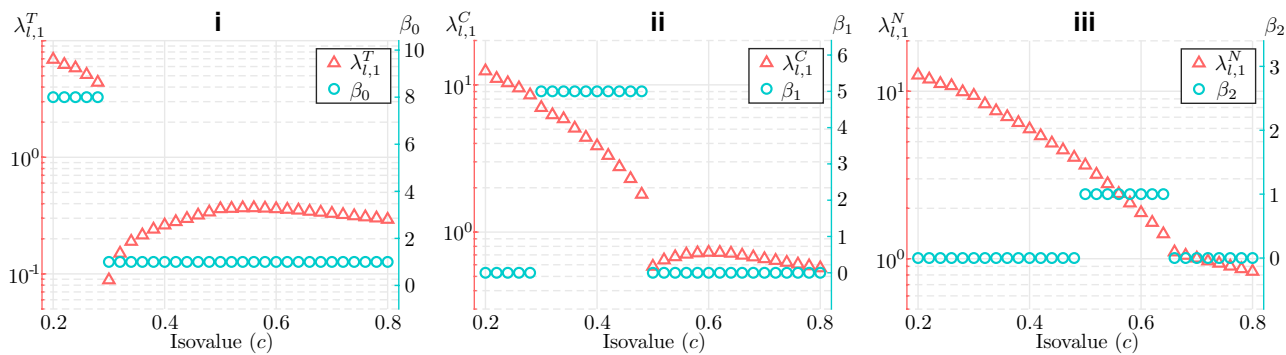
\begin{figure}
	\centering
	\includegraphics[height=1.82in]{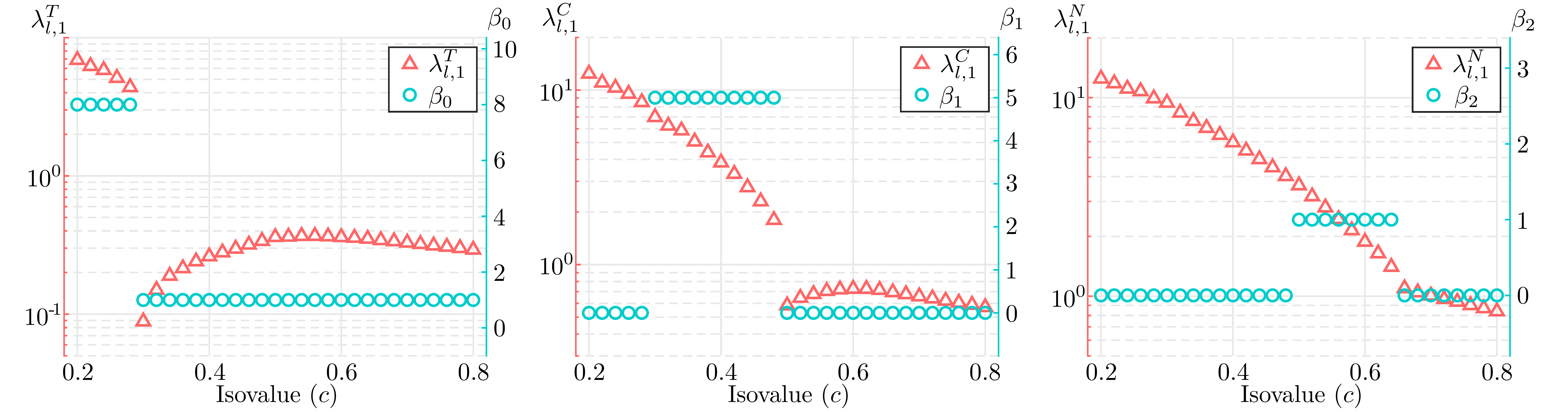}
	\begin{picture}(0,0)
	\put(-160, 135){\textbf{i}}
	\put(   5, 135){\textbf{ii}}
	\put( 165, 135){\textbf{iii}}
	\end{picture}
	\caption{ Eigenvalues and Betti numbers  vs isovalue ($c$) of the eight-body system with $\eta = 1.53$ and $\max(\rho)\approx1.1$. 
	{\bf i} shows the Fiedler values of the $T$ set and persistent Betti-0 numbers. 
	{\bf ii} shows the Fiedler values of the $C$ set and persistent Betti-1 numbers. 
	{\bf iii} illustrates the comparison of $\lambda_{l,1}^C$ and persistent $\beta_2$.
	}
	\label{fgr:eightAtomsEigBetti}
\end{figure}

We constructed the simple eight-body system to analyze the behavior of Hodge Laplacian spectra with an evolving cavity in the filtration. In this system, not only multiple connected components and multiple tunnels are involved, but a cavity also appears after the isovalue reaches a certain level before disappearing eventually. Thus, the dimension-2  Betti number $\beta_2$, which measures the number of cavities, changes during this process.

As shown in Fig.~\ref{fgr:eightAtoms}, the eight symmetric components start as blobs around eight vertices of a cube. Then they expand as the isovalue increases until they touch each other and form 6 rings, one for each face of the cube. At this point, persistent $\beta_0$ drops from 8 to 1, when  persistent $\beta_1$ increases from 0 to 5 (as five of the six tunnels are independent homology generators). As the level set value increases to the point that the tunnels are filled, persistent $\beta_1$ drops back to 0, but persistent $\beta_2$ increases to $1$ as a cavity formed inside the manifold. The cavity is filled up eventually, and persistent $\beta_2$ drops back to 0.

In Fig.~\ref{fgr:eightAtomsEigBetti}, the Fiedler values as functions of isovalue are shown in Figs.~\ref{fgr:eightAtomsEigBetti} {\bf i} and {\bf ii}, which exhibit similar behaviors as in the first two examples.
As in the previous example, the comparison between Figs.~\ref{fgr:eightAtomsEigBetti} {\bf i} and {\bf ii} shows that at $c=0.3$ the spectral function $\lambda^T_{l,1}$ starts to decay when $\lambda^C_{l,1}$ is discontinuous. Different from the previous examples, the smallest eigenvalues in {\bf iii} is no longer differentiable as persistent $\beta_2$ changes from one to zero near  isovalue 0.5. Fig.~\ref{fgr:eightAtomsEigBetti} {\bf iii} also indicates that at the isovalue where $\lambda^N_{l,1}$ is non-differentiable, $\lambda^C_{l,1}$ starts to decrease. Moreover, the simultaneous topological changes, the disappearance of tunnels and the appearance of the cavity, can be observed in $\lambda^C_{l,1}$. The disappearance of the cavity can be observed from $\lambda^N_{l,1}$.
From these preliminary results of the evolutionary de Rham-Hodge method, one may observe that the singular values in different spectra taken as functions of the isovalue $c$ not only illustrate the changes of topological features of different dimensions throughout the evolution of the manifold but also reveal the geometric features in different dimensions.
Therefore, empirically, the importance of low frequencies rather than the multiplicity of the zeroth frequency can already be observed in these simplistic constructions for features of different dimensionality. In the following, we demonstrate similar characteristics of spectral functions in two   molecular systems.

\subsection{Benzene molecule}

Benzene ($\text{C}_6\text{H}_6$) is a small organic chemical compound which consists of six carbon atoms in a planar hexagon ring and six hydrogen atoms each connected with one carbon atom. In this system, atoms have
different van der Waals radii, one for carbon and another for hydrogen. The carbon atoms are closer to each other than the hydrogen atoms and form the benzene ring. Thus, benzene is a perfectly simple yet realistic example to illustrate the evolutionary de Rahm-Hodge method. With the benzene data, we use $\eta=0.45$ to generate evolving manifolds.

\begin{figure}
	\centering
	\includegraphics[height=3.1in]{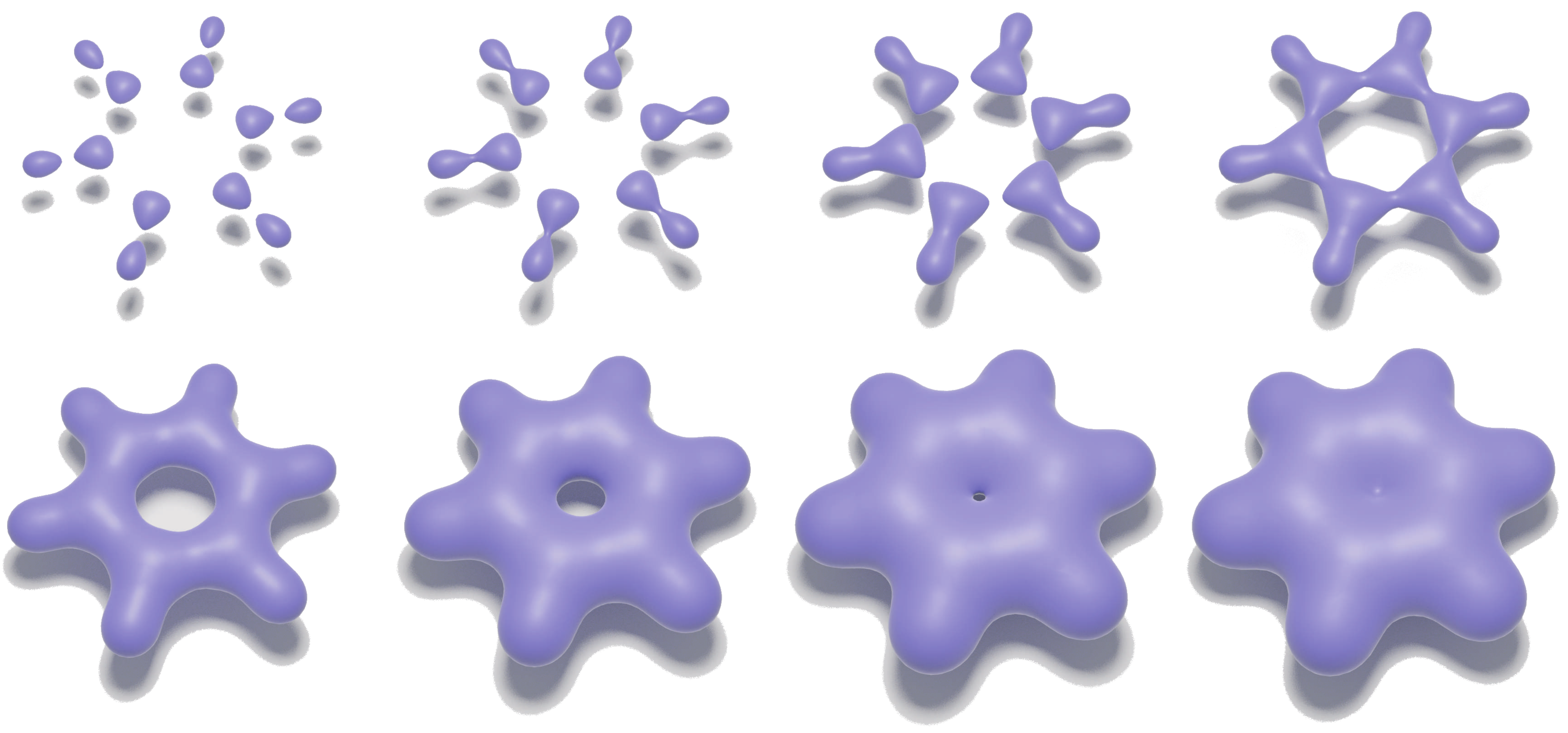}
	\begin{picture}(0,0)
	\put( -480, 210){\textbf{a}}
	\put( -360, 210){\textbf{b}}
	\put( -240, 210){\textbf{c}}
	\put( -120, 210){\textbf{d}}
	\put( -480, 110){\textbf{e}}
	\put( -360, 110){\textbf{f}}
	\put( -240, 110){\textbf{g}}
	\put( -120, 110){\textbf{h}}
	\end{picture}
	\caption{ Manifold evolution of benzene with $\eta=0.45\times r_{\rm vdw}$.
		{\bf a} through {\bf h} are snapshots from the start to the end.
		{\bf a} and {\bf b} show the transition of the persistent Betti-0 number from 12 to 6. 
		{\bf c} and {\bf d} show the formation of a ring;
		The Betti-0 number changes from 6 to 1 and remains at one to the end, whereas
		the Betti-1 number changes from zero to one.
		{\bf d}, {\bf e}, {\bf f} and {\bf g} illustrate the deformation of the hexagonal tunnel to a round tunnel.
		From {\bf g} to {\bf h}, the ring disappears and the Betti-1 number changes from 1 back to 0.
	}
	\label{fgr:benzene45}
\end{figure}
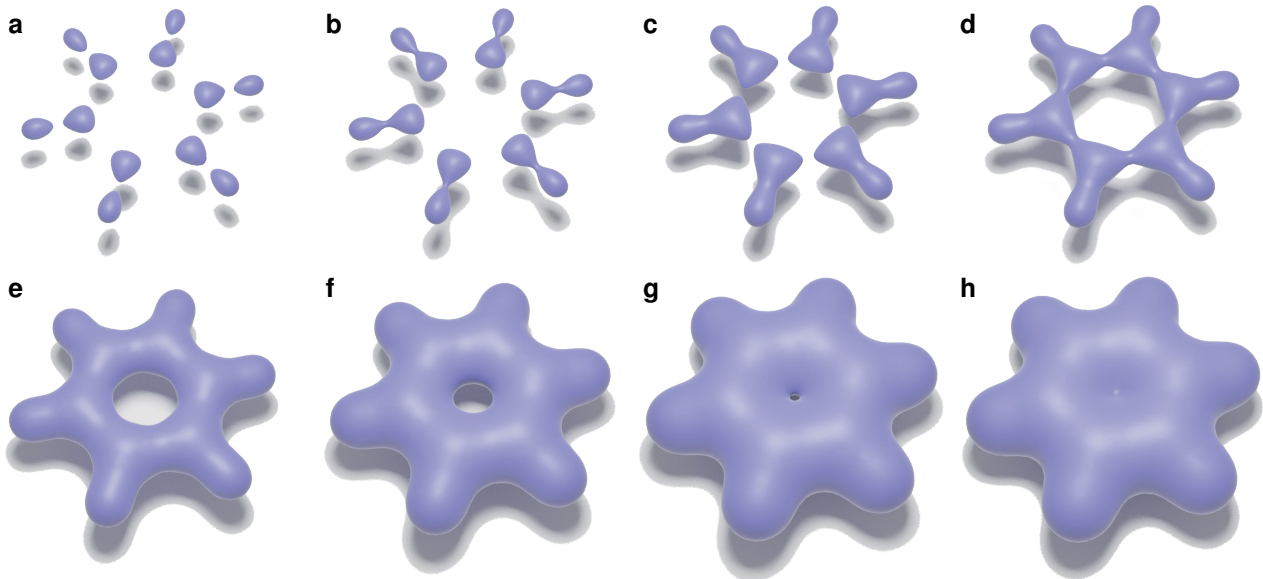

The first evolving manifold of benzene is generated at $\eta=0.45$. In the beginning, there are 12 components, with each smooth component center around one atom location as shown in Fig.~\ref{fgr:benzene45} {\bf a}. The van der Waals radius of carbon atoms is larger than that of hydrogen atoms, so the components associated with the carbon atoms are larger. From Fig.~\ref{fgr:benzene45} {\bf b} to Fig.~\ref{fgr:benzene45} {\bf c}, the originally separated components of the atoms start to connect pairwise, with a narrow tube formed between each hydrogen to its bonded carbon and thus, the persistent Betti-0 number is reduced to 6. The behavior of the manifold is similar to essentially six copies of our first example, the two-body system, until the six components of Fig.~\ref{fgr:benzene45} {\bf c} start to form a hexagonal ring, as shown in Fig.~\ref{fgr:benzene45} {\bf d}. At this point, there are six narrow tubes, one for each bond between two adjacent carbon atom pairs. As the density function continues to expand, the hexagonal ring evolves into a round cycle around a tunnel with a shrinking diameter. As the diameter of the tunnel reduces to zero at some parameter value between those of Fig.~\ref{fgr:benzene45} {\bf g} and Fig.~\ref{fgr:benzene45} {\bf h}, the noncontractible cycle disappears. During this topological change, the tiny cycle in the middle of the manifold in Fig.~\ref{fgr:benzene45} {\bf g} is filled up to form two concave surface patches in the middle of the manifold in Fig.~\ref{fgr:benzene45} {\bf h}.  The final topology of this system remains as a single component with a volume larger than that of Fig.~\ref{fgr:benzene45} {\bf h}.
 
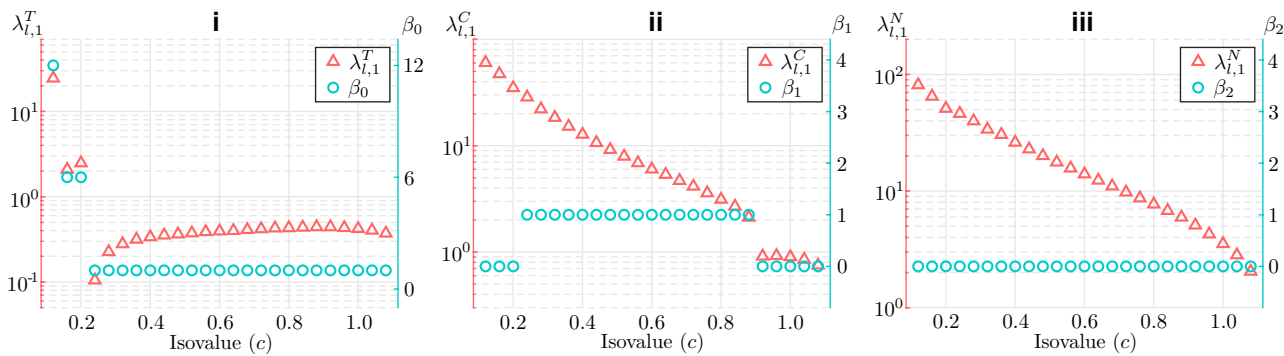
\begin{figure}
\centering
\includegraphics[height=1.82in]{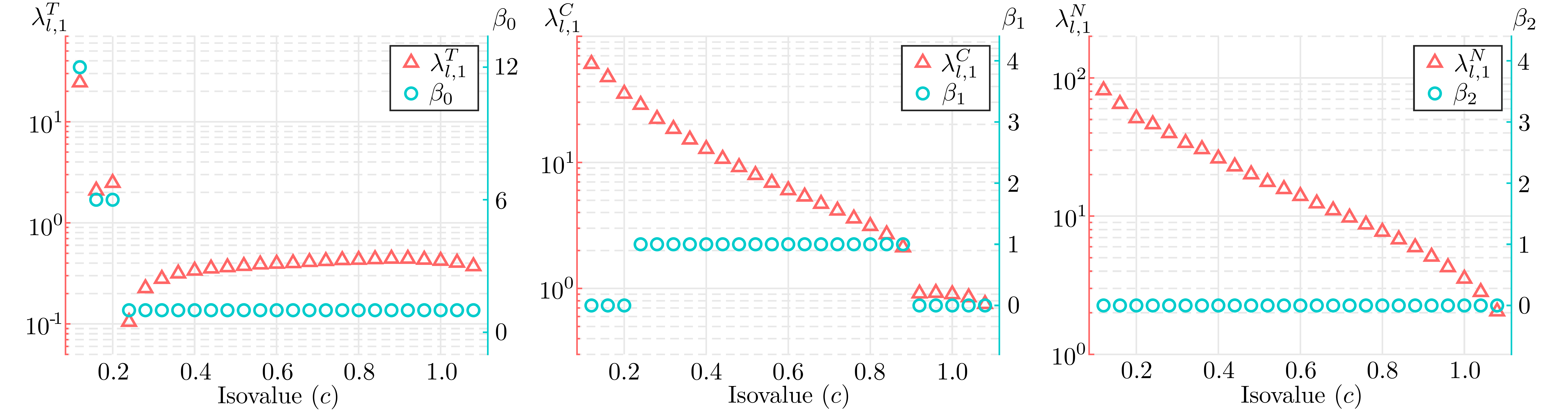}
\begin{picture}(0,0)
\put(-160, 135){\textbf{i}}
\put(   5, 135){\textbf{ii}}
\put( 165, 135){\textbf{iii}}
\end{picture}
\caption{ Eigenvalues and Betti numbers  vs isovalue ($c$) of the benzene system with $\eta = 0.45$ and $\max(\rho)\approx1.1$.
{\bf i} shows the smallest eigenvalues of the $T$ set. 
The drops at $c=0.12$ correspond to snapshots in Figs.~\ref{fgr:benzene45} {\bf a} and {\bf b}. 
The drops at $c=0.22$ correspond to snapshots in Figs.~\ref{fgr:benzene45} {\bf c} and {\bf d}.
{\bf ii} shows the smallest eigenvalue of the $N$ set.
The drops at $c=0.9$ correspond to snapshots in Figs.~\ref{fgr:benzene45} {\bf g} and {\bf h}.
{\bf iii} shows the smallest eigenvalues of the $C$ set.
}
\label{fgr:benzene45EigBetti}
\end{figure}

 Fig.~\ref{fgr:benzene45EigBetti} shows the Fiedler values of the $T$, $N$ and $C$ sets and their relations with the persistent Betti numbers when seen as a function of varying isovalues. First, for the $T$ set, $\lambda^T_{l,1}$ has two jumps at $c=0.12$ and $c=0.22$, which divide the $\lambda^T_{l,1}$ to three curve segments. Both discontinuities correspond to the decreases of the persistent Betti 0, from twelve to six, and then to one. As shown in Figs.~\ref{fgr:benzene45EigBetti} {\bf i}, $\lambda_{l,1}^T$ cannot only tell the topological changes but also give some additional information of a continuous portion of the evolution.  After $c=0.22$, $\lambda_{l,1}^T$ increases first and reaches its maximum at $c=0.9$ when the ring just disappears, at which point the structure (for tangential gradients) starts to grow softer as an expanding blob instead of a thicker ring. Fig.~\ref{fgr:benzene45EigBetti} {\bf ii} presents the jump of $\lambda_{l,1}^C$, which is correlated to the disappearance of the hole as indicated by the change of Betti-1 number from one to zero.
After the jump, $\lambda_{l,1}^C$ also increases slightly first and decays in the end. There is no cavity involved, so the spectral function shows a steady progression for the $C$ set as in our four-body example. One difference from that example is the finer grid used in the calculation, in order to handle the initial small components for the hydrogen atoms.

\subsection{Buckminsterfullerene}
The buckyball ($\text{C}_{60}$) has a beautiful structure composed of sixty carbon atoms.
It has twenty hexagons and twelve pentagons that resemble the pattern on a soccer ball, which has a rich structure with both geometric symmetries and  topology features. With our continuous density function, at certain values of $\eta$, the manifold evolution covers all the possible values of the persistent Betti-1 number allowed by the symmetry.
However, it is difficult to cover all the topological space for a density function associated with a single kernel size $\eta$. Thus we propose to use a multiscale (with a few different kernel sizes) analysis of the manifold evolution.
By using different $\eta$'s to capture different sets of snapshots for the evolving manifolds, we can compare the spectra across different kernel sizes $\eta$ as well as different control parameters $c$. We use the buckyball as an example for the multiscale analysis of  manifold evolution, and demonstrate how the spectra provide information on the evolution of their topological spaces and geometric features.

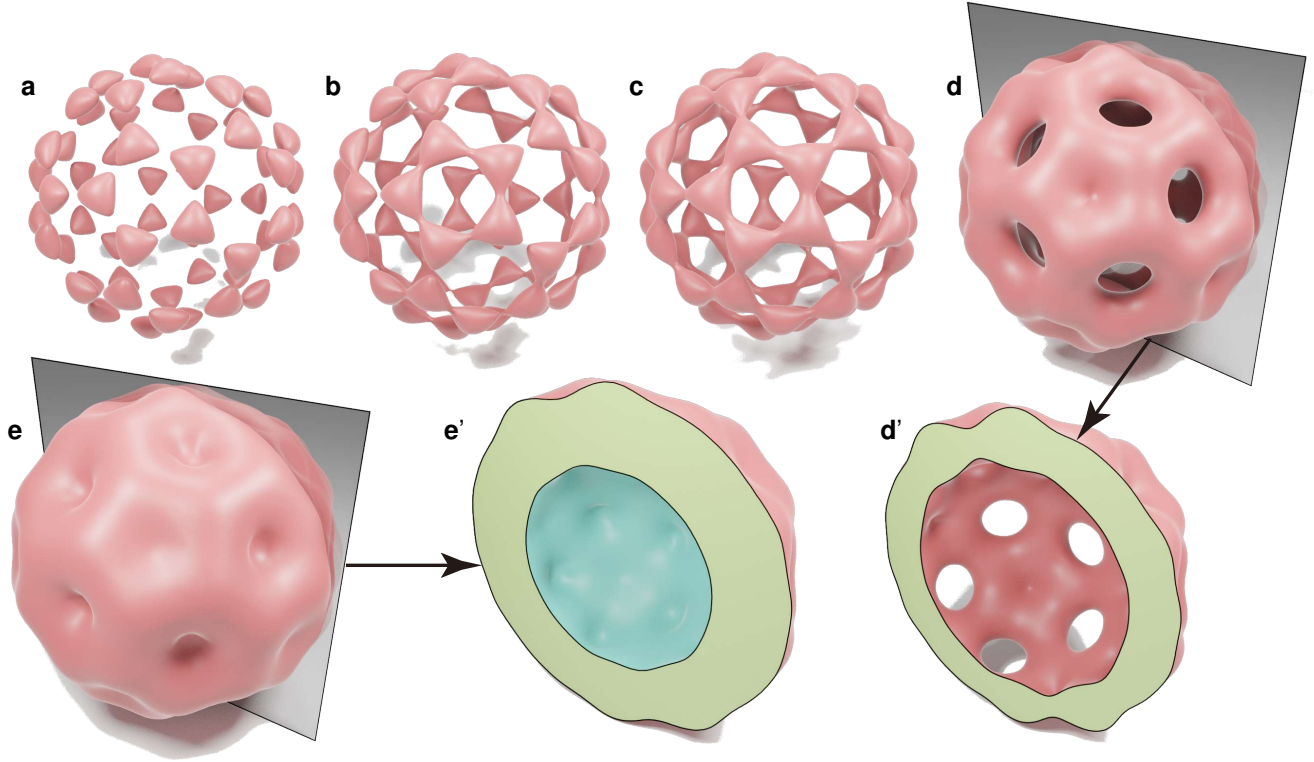
\begin{figure}
\centering
\includegraphics[height=4in]{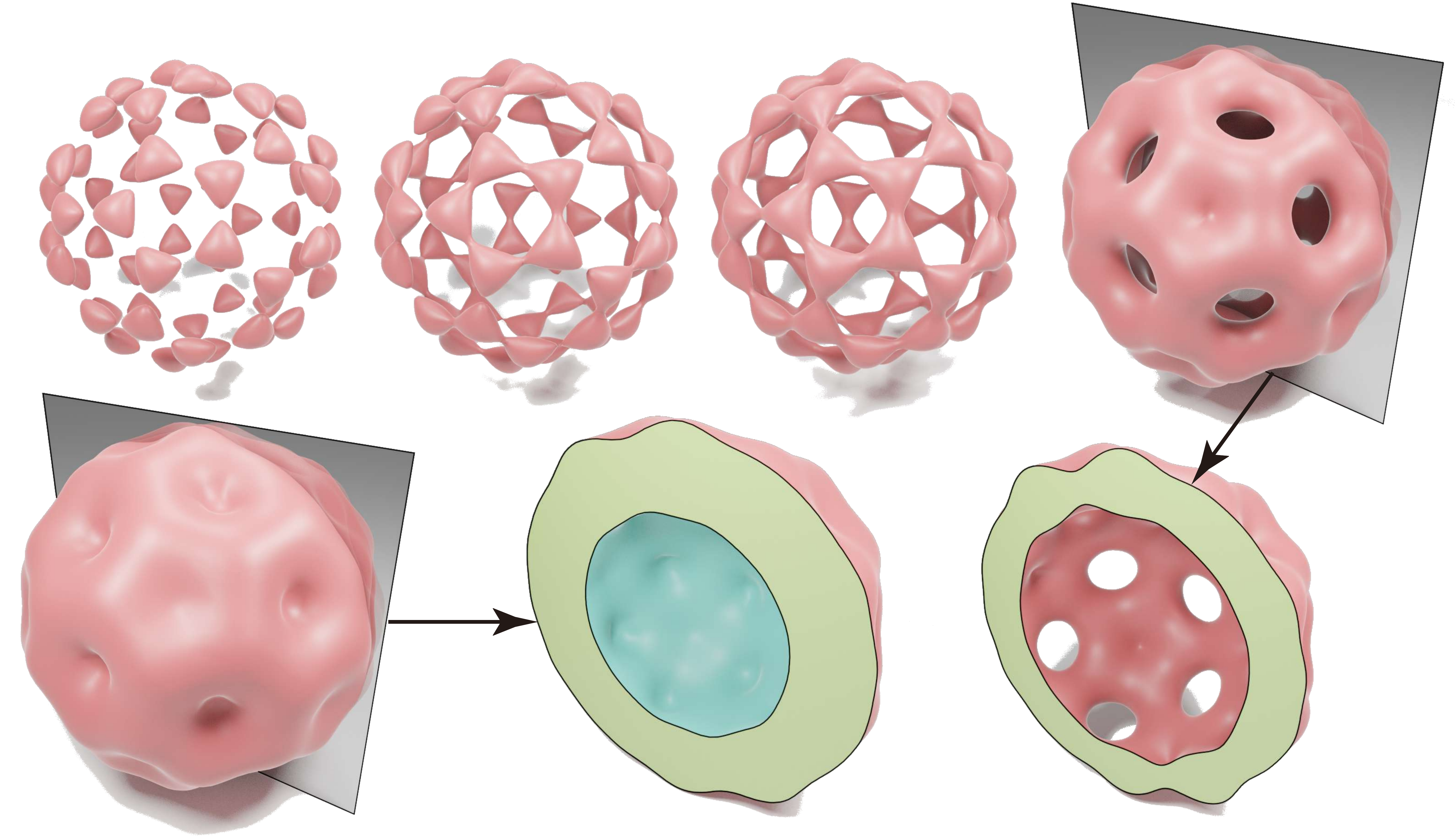}
\begin{picture}(0,0)
\put( -235, 265){\textbf{a}}
\put( -120, 265){\textbf{b}}
\put( -5, 265){\textbf{c}}
\put( 115, 265){\textbf{d}}
\put( 90, 135){\textbf{d}'}
\put( -240, 135){\textbf{e}}
\put( -75, 135){\textbf{e}'}
\end{picture}
\caption{Illustration of fullurene ($\text{C}_{60}$) manifold evolution with $\eta=0.5\times r_{\rm vdw}$.
{\bf a} presents sixty components around  carbon atom positions.
{\bf a} and {\bf b} show that the components connect if they share a pentagonal hole, and persistent $\beta_0$ changes from 60 to 12 and persistent $\beta_1$ changes from 0 to 12.
{\bf c} shows the hexagonal holes are formed, resulting in the change of persistent $\beta_0$ to 1 and persistent $\beta_1$ to 31. (There are 32 rings, but only 31 are independent in terms of homology.)
{\bf c} and {\bf d} show that the 12 pentagonal rings disappear and the persistent Betti-1 number drops from 31 to 19.
{\bf d} and {\bf e} show that the 20 hexagonal rings disappear and a cavity forms inside, so that persistent $\beta_1$ drops to $0$ and persistent $\beta_2$ increases to 1.
The vertical plan cuts the manifolds that gives an illustration of cavity in {\bf d'} and {\bf e'}.
}
\label{fgr:fullerene50}
\end{figure}
\begin{figure}
	\centering
	\includegraphics[height=1.82in]{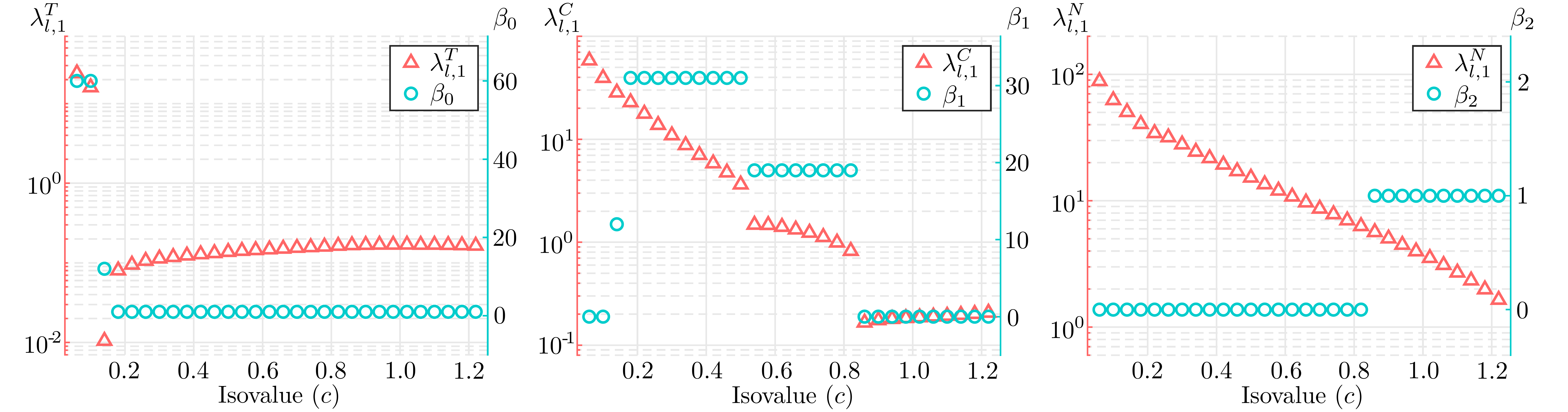}
	\begin{picture}(0,0)
	\put(-160, 135){\textbf{i}}
	\put(   5, 135){\textbf{ii}}
	\put( 165, 135){\textbf{iii}}
	\end{picture}
	\caption{Eigenvalues and Betti numbers  vs isovalue ($c$) of the fullurene ($\text{C}_{60}$) system with $\eta = 0.5\times r_{\rm vdw}$ and $\max(\rho)\approx1.3$.
		{\bf i} gives the Fiedler values of the $T$ set and persistent $\beta_0$.
		{\bf ii} presents the comparison of $\lambda^C_{l,1}$ and persistent $\beta_1$.
		{\bf iii} shows the Fiedler values of the $N$ set and persistent $\beta_2$.
	}
	\label{fgr:fullerene50EigBetti}
\end{figure}

For kernel scaling parameter $\eta=0.5\times r_{\rm vdw}$, the manifold evolution starts with 60 components as shown in Fig.~\ref{fgr:fullerene50} {\bf a}. The components start the expansion, each around the position of one carbon atom, and merge into larger connected components if they share a common pentagon in the skeleton structure as shown in Fig.~\ref{fgr:fullerene50} {\bf b}. This leads to the changes in persistent $\beta_0$ (from 60 to 12) and persistent $\beta_1$ (from 0 to 12). Fig.~\ref{fgr:fullerene50} {\bf c} shows the snapshot right after the appearance of twenty hexagonal holes. Next, each hole starts to shrink. As each pentagonal hole has a smaller size than that of a hexagonal hole, we observe in Fig.~\ref{fgr:fullerene50} {\bf c} to Fig.~\ref{fgr:fullerene50} {\bf d}, the pentagonal holes disappear before the hexagonal holes also disappear. Simultaneous to the disappearance of hexagons, a cavity is created. In Fig.~\ref{fgr:fullerene50} {\bf e} after the formation of the cavity, both the outer surface and the inner surface contain numerous regions of concavity and gradually, the shape evolves to resemble a slightly dented thick spherical shell.

For analysis of this evolution, Fig.~\ref{fgr:fullerene50EigBetti} illustrates the eigenvalues and Betti numbers versus the isolvaue $c$. Fig.~\ref{fgr:fullerene50EigBetti} {\bf i} gives the Fiedler values (smallest eigenvalue) of the $T$ set and $\beta_0$. This Betti number has two drops, from 60 to 12, and then to 1. Within each interval of isovalues with the same persistent Betti number, $\lambda_{l,1}^T$ is changing smoothly as expected from our discussion on homeomorphic shapes with a slowly evolving metric. Fig.~\ref{fgr:fullerene50EigBetti} {\bf ii} presents the information that the Fiedler values of the $C$ set can offer. For the interval, $c\in [0.16,0.5]$,  persistent $\beta_1$ remains at $31$, and the continuous decrease in $\lambda^C_{l,1}$ shows that the geometric structure is ``softer'' for the curl fields as the handles grow thicker. Similarly, for intervals within which persistent $\beta_1$ equals to 19 or 1, $\lambda^C_{l,1}$ is a smooth function within each interval but is discontinuous at the boundary of these intervals where the topology transitions.
The Fiedler values of the $N$ set are given in Fig.~\ref{fgr:fullerene50EigBetti} {\bf iii}, which, although mostly smooth, also has changed in slope at isovalues associated with changes in connected components and tunnels.
As the examples become more complex, the spectral functions also exhibit richer structure, with the advantage of indicating both topological persistence and geometric progression.

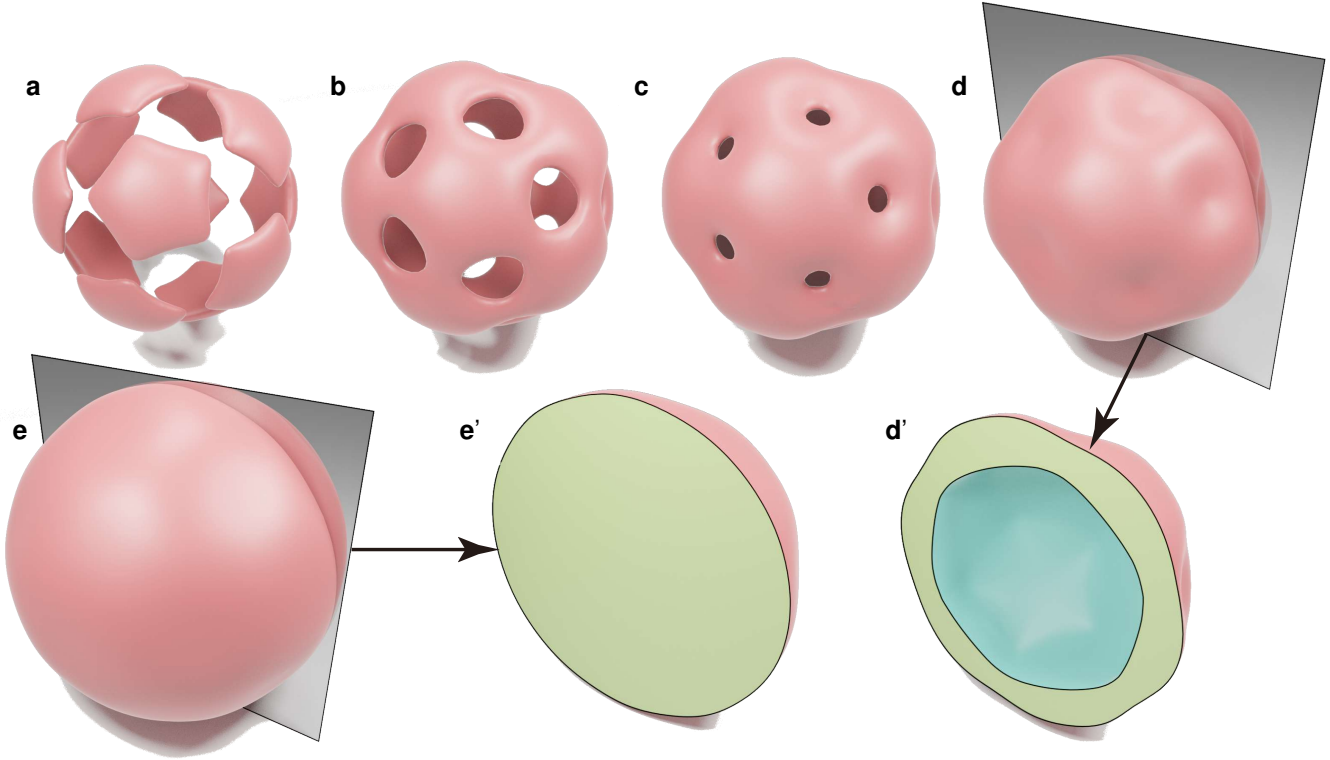
\begin{figure}
	\centering
	\includegraphics[height=4in]{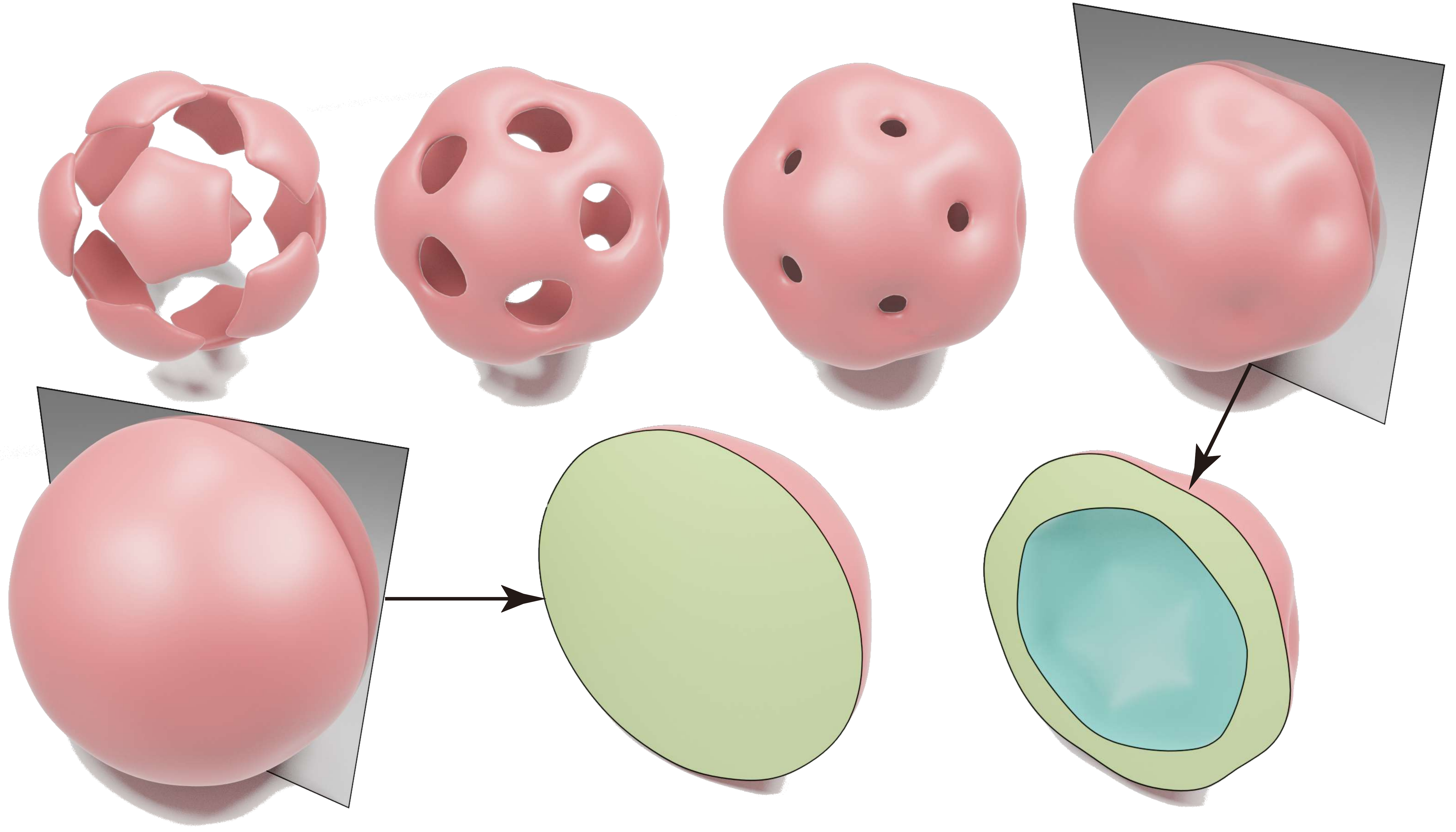}
	\begin{picture}(0,0)
	\put( -235, 265){\textbf{a}}
	\put( -120, 265){\textbf{b}}
	\put( -5, 265){\textbf{c}}
	\put( 115, 265){\textbf{d}}
	\put( 90, 135){\textbf{d}'}
	\put( -240, 135){\textbf{e}}
	\put( -71, 135){\textbf{e}'}
	\end{picture}
	\caption{Illustration of fullurene ($\text{C}_{60}$) manifold evolution with $\eta=0.8\times r_{\rm vdw}$.
		{\bf a} shows 12 initial solid pentagonal components.
		{\bf b} and {\bf c} show the formation and contraction process of the 20 rings.
		{\bf d} is the snapshot right after the formation of the cavity.
		{\bf e} shows the final stage as a solid ball of this example.
	}
	\label{fgr:fullerene80}
\end{figure}
\begin{figure}
	\centering
	\includegraphics[height=1.82in]{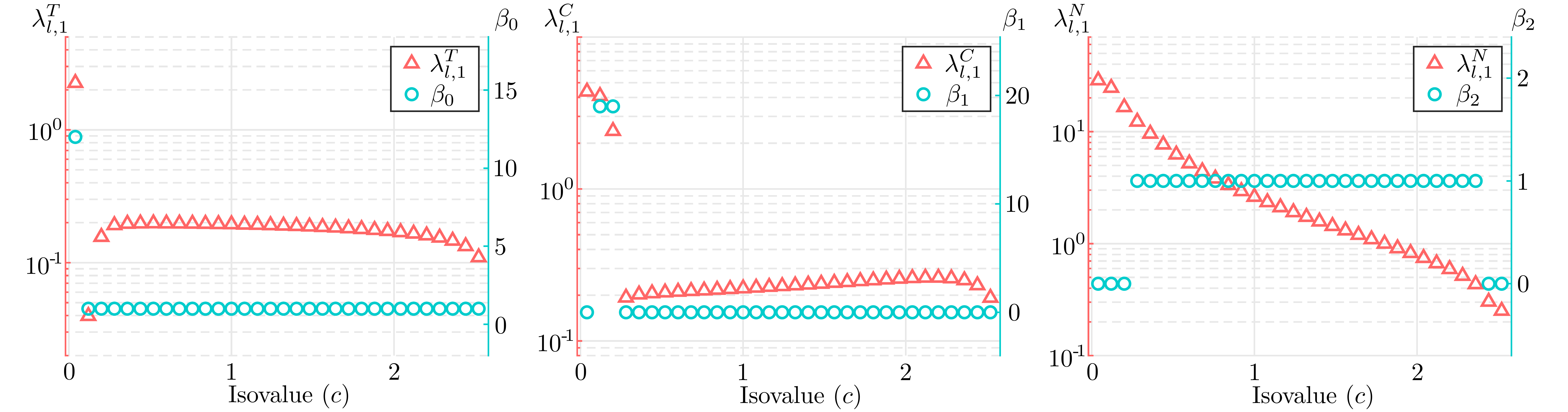}
	\begin{picture}(0,0)
	\put(-160, 135){\textbf{i}}
	\put(   5, 135){\textbf{ii}}
	\put( 165, 135){\textbf{iii}}
	\end{picture}
	\caption{ Eigenvalues and Betti numbers  vs isovalue ($c$) of the fullurene ($\text{C}_{60}$) system with $\eta = 0.8\times r_{\rm vdw}$; $\max{\rho}\approx2.5$.
	{\bf i} gives the Fiedler values of the $T$ set and persistent $\beta_0$.
	{\bf ii} presents the comparison of $\lambda^C_{l,1}$ and  persistent $\beta_1$.
	{\bf iii} shows the Fiedler values of the $N$ set and persistent $\beta_2$.
	}
	\label{fgr:fullerene80EigBetti}
\end{figure}

For large and dense point sets as in this fullerene, the shape of the manifold evolution is heavily influenced by the kernel size $\eta$. To show the importance of multiscale analysis, we create a second evolution with $\eta=0.8\times r_{\rm vdw}$ and generate the snapshots in Fig.~\ref{fgr:fullerene80}. For the initial isovalue, as seen in Fig.~\ref{fgr:fullerene80} {\bf a}, the manifold consists of twelve pentagonal components. Unlike the evolution with $\eta=0.5\times r_{\rm vdw}$, which contains pentagonal holes alongside hexagonal holes, here the pentagonal components are already with the holes filled before the hexagonal holes are even formed. Thus, the two evolutions cannot find a homeomorphism between their stages even if any isovalues are allowed, which implies that they can reveal different information regarding the system. As the components connect, twenty rings show up as in Figs.~\ref{fgr:fullerene80} {\bf b} and \ref{fgr:fullerene80} {\bf c}, with decreasing diameters for increasing isovalues. Once the cavity is formed, the large inner surface shown in Fig.~\ref{fgr:fullerene80} {\bf d} starts to contract, and the manifold ends up as a solid ball in Fig.~\ref{fgr:fullerene80} {\bf e}. As for the spectral functions, Fig.~\ref{fgr:fullerene80EigBetti} shows three plots of the Fiedler values of the $T$, $C$ and $N$ sets and the persistent Betti numbers against the isovalues, respectively. Since the components connect right after first two snapshots, Fig.~\ref{fgr:fullerene80EigBetti} {\bf i} shows the drop of $\lambda_{l,1}^T$  in the third snapshot as persistent $\beta_1$ changes from 12 to 1. The Fiedler values $\lambda_{l,1}^T$ then increases before starting to decrease when persistent $\beta_1$ drops to 0 when the system can be seen as a shell growing softer with thicker membrane instead of a structure growing stiffer with thicker supporting handles.  Similarly, there are only a few snapshots for the evolving manifold to have rings as they are quickly filled up. In Fig.~\ref{fgr:fullerene80EigBetti} {\bf ii}, the Fiedler values $\lambda_{l,1}^C$ already decreases quickly before plunging to a small number at the point when holes disappear. During the period of the inner surface contracting and outer surface expanding, $\lambda_{l,1}^C$ increases first as the structure grows stiffer for curl fields, and then grows softer eventually near the very end of the manifold evolution. In the last plot of Fig.~\ref{fgr:fullerene80EigBetti}, $\lambda_{l,1}^N$ slightly increases at beginning and then decreases smoothly. The disappearance of the cavity is captured at the end of snapshots, thus there is a non-differentiable point at end of this spectral function. We see in this evolution again, that the progression of the manifold evolution can be observed in the spectral functions as well as the topological transitions.

\section{Conclusion} \label{conclusion}
While persistent homology has had tremendous success in data science and machine learning via a multiscale analysis, it does not capture geometric progression when there are no topological changes. In contrast, although de Rham-Hodge theory provides a simultaneous geometric and topological analysis, it lacks multiscale information. We introduce an evolutionary de Rham-Hodge method to offer a unified multiscale geometric and topological representation of data. The evolutionary de Rham-Hodge method is applied to analyze the topological and geometric characteristics through the evolution of manifolds which are a family of 3D multiscale shapes constructed from an evolutionary filtration process. In addition to exactly the topological persistence that would be obtained from persistent homology, the analysis of the evolutionary spectra of Hodge Laplacian operators portrays geometric progression. Specifically,  appropriate treatments of the  Hodge Laplacian boundary conditions gives rise to three unique sets of singular spectra associated with the tangential gradient eigen field ($T$), the curl eigen field ($C$), and the tangential divergent eigen field ($N$). The multiplicities of the  zero eigenvalues  corresponding to  the $T$, $C$, and  $N$ sets of spectra are exactly  the persistent Betti-0 ($\beta_0$),  Betti-1 ($\beta_1$), and Betti-2 ($\beta_2$) numbers one would obtain from persistent homology. Using discrete exterior calculus in close manifolds or compact manifolds with boundary, we show that investigating the first non-zero eigenvalues, i.e., Fiedler values, of the $T$, $C$, and  $N$ sets of evolutionary spectra unveil both the persistence for topological features and the geometric progression for the shape analysis. For a proof-of-concept analysis, the evolutionary de Rham-Hodge method is applied to a few benchmark examples, including the two-body system, four-body system, eight-body system,  benzene (C$_{6}$H$_{6}$), and buckminsterfullerene (C$_{60}$). Extensive numerical experiments demonstrate that the present evolutionary de Rham-Hodge method captures the multiscale geometric progression and topological persistence of data.

In the present proof-of-concept analysis, only the first non-zero eigenvalues are presented. However, in practical applications, both eigenfunctions and high-order eigenvalues are needed in de Rham-Hodge  modeling and analysis as shown in our recent work \cite{zhao2019rham}. The proposed evolutionary de Rham-Hodge method provides a solid foundation for a wide variety of applications,  including shape analysis, image processing, computer vision, pattern recognition, computer aided design, network analysis,  computational biology, and drug design. Since the evolutionary de Rham-Hodge method can reveal both topological persistence and geometric progression, it will offer a powerful multiscale representation of data for machine learning, including deep learning.   

Finally, the present evolutionary de Rham-Hodge method opens new opportunities in further theoretical developments in differential geometry,  such as the introduction of multiscale analysis to Riemannian connection, tensor bundle, index theory,  and K theory.

\section{Acknowledgments}
This work was supported in part by NSF Grants DMS1721024, DMS1761320, IIS1900473, NIH grants GM126189 and GM129004,  Bristol-Myers Squibb,  and Pfizer. GWW thanks Vidit Nanda for useful discussion.


\bibliographystyle{alpha}

\end{document}